\documentstyle [12pt]{article}
\def\picture#1by#2(#3){
\vbox to #2 {
  \hrule width #1 height 0pt depth 0pt \vfill \special{picture #3}}
}

\def\scaledpicture#1by#2(#3scaled#4){{
\dimen0=#1  \dimen1=#2
\divide\dimen0 by 1000 \multiply\dimen0 by #4
\divide\dimen1 by 1000 \multiply\dimen1 by #4
\picture \dimen0 by \dimen1 (#3 scaled #4)}}
\def\dfigure#1by#2(#3scaled#4offset#5:#6)
  {\medskip
   \vglue 2mm minus 2mm
   $$
     \hbox{
       \hglue#5
       {\scaledpicture #1 by #2 (#3 scaled #4)}
     }
   $$
   \par\goodbreak
   \vglue 2mm minus 2mm
   \medskip}

\def\qmod#1#2{{\hbox{}^{\displaystyle{#1}}}\!\big/\!\hbox{}_{
\displaystyle{#2}}}

\newfam\msbfam
\font\twelmsb=msbm10 at 12pt

\font\sevenmsb=msbm10 at 7pt
\font\fivemsb=msbm10 at 5pt

\textfont\msbfam=\twelmsb
\scriptfont\msbfam=\sevenmsb
\scriptscriptfont\msbfam=\fivemsb

\def\Bbb{\fam\msbfam\twelmsb}

\def\C{{\Bbb C}}

\def\G{{\Bbb G}}
\def\H{{\Bbb H}}

\def\N{{\Bbb N}}

\def\P{{\Bbb P}}

\def\R{{\Bbb R}}
\def\Z{{\Bbb Z}}

\def\qed {\hfill\vrule height6pt width6pt depth0pt \bigskip}

\def\map{\longrightarrow}
\def\textmap#1{\mathop{\vbox{\ialign{
                                ##\crcr
    ${\scriptstyle\hfil\;\;#1\;\;\hfil}$\crcr
    \noalign{\kern-1pt\nointerlineskip}
    \rightarrowfill\crcr}}\;}}

\def\textlmap#1{\mathop{\vbox{\ialign{
                                ##\crcr
    ${\scriptstyle\hfil\;\;#1\;\;\hfil}$\crcr
    \noalign{\kern-1pt\nointerlineskip}
    \leftarrowfill\crcr}}\;}}

\newfam\meuffam
\font\twelvmeuf=eufm10 at 12 pt
\font\tenmeuf=eufm10
\font\sevenmeuf=eufm7

\textfont\meuffam=\twelvmeuf
\scriptfont\meuffam=\tenmeuf
\scriptscriptfont\meuffam=\sevenmeuf

\def\germ{\fam\meuffam\tenmeuf}

\def\cg{{\germ c}}

\def\g{{\germ g}}

\def\mg{{\germ m}}

\begin{document}
\def\Pr{{\rm Pr}}
\def\tr{{\rm Tr}}
\def\End{{\rm End}}
\def\Spin{{\rm Spin}}
\def\U{{\rm U}}
\def\SU{{\rm SU}}
\def\SO{{\rm SO}}
\def\PU{{\rm PU}}
\def\GL{{\rm GL}}
\def\spin{{\rm spin}}
\def\u{{\rm u}}
\def\su{{\rm su}}
\def\so{{\rm so}}
\def\ub{\underbar}
\def\pu{{\rm pu}}
\def\Pic{{\rm Pic}}
\def\Iso{{\rm Iso}}
\def\NS{{\rm NS}}
\def\deg{{\rm deg}}
\def\Hom{{\rm Hom}}
\def\h{{\germ h}}
\def\Herm{{\rm Herm}}
\def\Vol{{\rm Vol}}
\def\pf{{\bf Proof: }}
\def\id{{\rm id}}
\def\i{{\germ i}}
\def\im{{\rm im}}
\def\rk{{\rm rk}}
\def\ad{{\rm ad}}
\def\h{{\bf H}}
\def\coker{{\rm coker}}
\def\dv{\bar\partial}
\def\Ad{{\rm Ad}}
\def\RSU{\R SU}
\def\ad{{\rm ad}}
\def\dva{\bar\partial_A}
\def\da{\partial_A}
\def\p{\partial\bar\partial}
\def\sp{\Sigma^{+}}
\def\sm{\Sigma^{-}}
\def\spm{\Sigma^{\pm}}
\def\smp{\Sigma^{\mp}}
\def\oo{{\scriptstyle{\cal O}}}
\def\ooo{{\scriptscriptstyle{\cal O}}}
\def\sw{Seiberg-Witten }
\def\pa{\partial_A\bar\partial_A}
\def\Dr{{\raisebox{0.15ex}{$\not$}}{\hskip -1pt {D}}}
\def\gr{{\scriptscriptstyle|}\hskip -4pt{\g}}
\def\subsetint{{\  {\subset}\hskip -2.45mm{\raisebox{.28ex}
{$\scriptscriptstyle\subset$}}\ }}
\def\nr{\parallel}
\newtheorem{sz}{Satz}[subsection]
\newtheorem{thry}[sz]{Theorem}
\newtheorem{pr}[sz]{Proposition}
\newtheorem{re}[sz]{Remark}
\newtheorem{co}[sz]{Corollary}
\newtheorem{dt}[sz]{Definition}
\newtheorem{lm}[sz]{Lemma}
\newtheorem{cl}[sz]{Claim}

\title{Recent Developments in Seiberg-Witten Theory and Complex Geometry}
\author{\\ Christian Okonek\thanks{Partially supported  by: AGE-Algebraic
Geometry in Europe, contract No ERBCHRXCT940557 (BBW 93.0187), and by  SNF, nr.
21-36111.92}
 \and  \\   Andrei Teleman$^*$  }

\date{  }
\maketitle

\setcounter{section}{-1}
\tableofcontents
\newpage
\hspace{2.5cm}{ \small We dedicate this paper to our wives Christiane and
Roxana for}\\
\hspace*{3cm} {\small  \
their invaluable help and support  during the past two years.}\\
\section{Introduction}
\vspace*{5mm}

About two years  ago, in October 1994, E. Witten revolutionized the theory
of 4-manifolds
by introducing the now famous \sw invariants [W]. These invariants are
defined by counting
gauge equivalence classes of solutions of the \sw monopole equations, a
system of
non-linear  PDE's which describe the absolute minima of a Yang-Mills-Higgs type
functional with an abelian gauge group.

In a very short period   of only a few weeks after Witten's seminal paper
became available,
several long-standing conjectures were solved, many new and totally
unexpected results
were found, and much simpler and more conceptional proofs of already established
theorems were given.

Among the most spectacular applications in this early period are the
solution  of the Thom
conjecture [KM], new results about Einstein metrics and Riemannian metrics
of positive
scalar curvature [L1] [L2], a proof of a  $10/8$ bound for intersection
forms of  $Spin$
manifolds [F], as well as  several results about the ${\cal
C}^{\infty}$-classification  of
algebraic surfaces [OT1], [OT2], [FM], [Bru].
The latter include Witten's proof of the ${\cal C}^{\infty}$-invariance of
the canonical
class  of a minimal surface of general type with $b_+\ne1$  up to sign, and
a  simple proof of
the Van de Ven conjecture by the authors
\footnote{Combining results in [L1], [L2] with
ideas from [OT2], P. Lupa\c scu recently obtained the optimal
characterization of complex
surfaces of K\"ahler type admitting Riemannian metrics of non-negative
scalar curvature
[Lu]}.

In two of the earliest papers on the subject, C. Taubes found a deep
connection between
\sw theory and symplectic geometry in dimension four: He first showed that
many aspects
of the new theory extend from the case of K\"ahler surfaces  to the more general
symplectic case [Ta1], and then he went  on to establish a beautiful
relation between \sw
invariants and Gromov-Witten invariants of symplectic 4-manifolds [Ta2].

A report on some papers of this first period can be found in [D].

Since the time this report was written, several new developments have taken
place:

The original \sw theory, as introduced in [W], has been refined and
extended to the case of
manifolds with $b_+=1$. The structure of the \sw invariants is  more
complicated in
this situation, since the invariants for manifolds with $b_+=1$ depend on a
chamber
structure. The general theory, including the complex-geometric
interpretation in the case
of K\"ahler surfaces, is now completely understood [OT6].

At present, three major
directions of research have emerged:\\
- Seiberg-Witten theory and symplectic geometry\\
- Non-abelian \sw theory and complex geometry\\
- Seiberg-Witten-Floer theory and contact structures

In this article, which had its origin in the notes for several lectures
which we gave in
Berkeley, Bucharest, Paris, Rome and Z\"urich during the past two years, we
concentrate
mainly on the second of these directions.

The reader will probably notice that the non-abelian theory  is a subject
of much higher
complexity than the original (abelian) \sw theory; the difference is
roughly comparable to
the difference between  Yang-Mills theory and Hodge theory. This complexity
accounts for
the length of the article. In rewriting our notes, we have tried to
describe the essential
constructions as simply as possible but without oversimplifying, and we
have made an
effort to explain the most important ideas and results carefully in a
non-technical way;
for proofs and technical details precise references are given.

We hope that this presentation of the material will motivate the reader,
and we believe
that our notes can serve as a comprehensive introduction to an interesting
new field of
research.

We have divided the article in three chapters. In chapter 1 we give a
concise but complete
exposition of the basics of abelian \sw theory in its most general form.
This includes the
definition of refined invariants for manifolds with $b_1\ne 0$, the
construction of
invariants for manifolds with $b_+=1$, and the universal wall crossing
formula in this
situation.

Using this formula in connection with vanishing and transversality results,
we calculate
the \sw  invariant for the simplest non-trivial example, the projective plane.

In chapter 2 we introduce non-abelian \sw theories for rather general
structure groups $G$.
After a careful exposition of $Spin^G$-structures and $G$-monopoles, and a short
description of some important properties of their moduli spaces, we explain
one of the main
results of the Habiliations\-schrift of the second author [T2], [T3]: the
fundamental
  Uhlenbeck type
compactification of the   moduli spaces of $PU(2)$-monopoles.

Chapter 3 deals with  complex-geometric aspects of \sw theory: We show that on
K\"ahler surfaces moduli spaces of $G$-monopoles, for unitary structure
groups $G$, admit
an interpretation as moduli spaces of purely holomorphic objects. This
result is a
Kobayashi-Hitchin type correspondence whose proof depends on a careful
analysis of the
relevant vortex equations. In the abelian case it identifies the moduli
spaces of twisted
\sw monopoles with certain Douady spaces of curves on the surface [OT1]. In the
non-abelian case we obtain an identification between moduli spaces of
$PU(2)$-monopoles
and moduli spaces of stable oriented pairs [OT5], [T2].

The relevant stability concept is new and makes sense on K\"ahler manifolds
of arbitrary
dimensions; it is induced by a natural moment map which is closely related
to the
projective vortex equation. We clarify the connection between this new
equation and the
parameter dependent vortex equations which had been studied in the
literature [Br]. In the
final section we  construct moduli spaces of stable oriented pairs on
projective varieties
of any dimension with GIT methods  [OST]. Our moduli spaces are projective
varieties which
come with a natural $\C^*$-action, and they play the role of master spaces
for stable
pairs.  We end our article  with the description of a very general
construction principle
which we call "coupling and reduction". This fundamental principle allows
to  reduce the
calculation of correlation functions associated with vector bundles to a
computation on the
space of reductions, which is essentially a moduli space of lower rank objects.

Applied to suitable master spaces on curves, our principle yields a
conceptional new proof
of the Verlinde formulas, and very likely also a proof of the
Vafa-Intriligator conjecture.
The gauge theoretic version of the same principle can be used to prove
Witten's conjecture,
and more generally, it will probably also lead to formulas expressing the
Donaldson invariants
of arbitrary 4-manifolds in terms of \sw invariants.

\section{Seiberg-Witten invariants}

\subsection{The monopole equations}

Let $(X,g)$ be a closed oriented Riemannian 4-manifold. We denote by
$\Lambda^p$ the bundle
of $p$-forms on $X$ and by $A^p:=A^0(\Lambda^p)$ the corresponding  space
of sections.
Recall that the Riemannian metric $g$ defines a Hodge operator $*:\Lambda^p\map
\Lambda^{4-p}$  with $*^2=(-1)^p$. Let
$\Lambda^2=\Lambda^2_+\oplus\Lambda^2_-$ be the corresponding eigenspace
decomposition.

 A \underbar{$Spin^c$}-\underbar{structure} on $(X,g)$ is a  triple
$\tau=(\Sigma^{\pm},\iota,\gamma)$ consisting of a pair of $U(2)$-vector bundles
$\Sigma^{\pm}$, a unitary isomorphism $\iota:\det\Sigma^+\map\det\Sigma^-$
and an
orientation-preserving linear isometry
$\gamma:\Lambda^1\map\RSU(\Sigma^+,\Sigma^-)$.
Here $\RSU(\Sigma^+,\Sigma^-)\subset \Hom_\C(\Sigma^+,\Sigma^-)$ is the
subbundle of real
multiples of (fibrewise) isometries of determinant 1.  The \underbar{spinor}
\underbar{bundles} $\spm$ of $\tau$ are -- up to isomorphism --  uniquely
determined by their
first Chern class $c:=\det\spm$, the Chern class of  the
$Spin^c(4)$-structure $\tau$.   This
class can be any integral lift of the second Stiefel-Whitney class $w_2(X)$
of $X$, and, given
$c$, we have
$$c_2(\spm)=\frac{1}{4}(c^2-3\sigma(X)\mp 2e(X)) \ .
$$
Here $\sigma(X)$ and $e(X)$ denote the signature and the Euler
characteristic of $X$.

The map $\gamma$ is called the \underbar{Clifford} \underbar{map} of the
$Spin^c$-structure
$\tau$.  We denote by $\Sigma$ the total spinor bundle
$\Sigma:=\sp\oplus\sm$, and we   use
the same symbol $\gamma$ also for the induced the map $\Lambda^1\map
su(\Sigma)$ given
by
$$  u\longmapsto\left(\matrix{0&-\gamma(u)^*\cr\gamma(u)&0}\right) \ .
$$
 Note that the Clifford identity
$$\gamma(u)\gamma(v)+\gamma(v)\gamma(u)=-2g(u,v)
$$
holds, and that the formula
$$\Gamma(u\wedge v):=\frac{1}{2}[\gamma(u),\gamma(v)]
$$
defines an embedding $\Gamma:\Lambda^2\map su(\Sigma)$ which maps
$\Lambda^2_{\pm}$ isometrically onto $su(\spm)\subset su(\Sigma)$.

The second cohomology group $H^2(X,\Z)$ acts on the set   of
 equivalence classes $\cg$ of  $Spin^c(4)$-structures on $(X,g)$ in a
natural way: Given a
representative $\tau=(\spm,\iota,\gamma)$ of $\cg$ and a Hermitian line
bundle $M$
representing a class $m\in H^2(X,\Z)$, the tensor product
$(\spm\otimes M,
\iota\otimes \id_{M^{\otimes 2}},\gamma\otimes\id_M)$ defines a
$Spin^c$-structure
$\tau_m$. Endowed with the   $ H^2(X,\Z)$-action given by
$(m,[\tau])\longmapsto  [\tau_m]$, the set of equivalence classes of
$Spin^c$-structures   on
$(X,g)$ becomes   a $H^2(X,\Z)$-torsor, which is independent of the metric
$g$  up to canonical
isomorphism [OT6]. We denote this $H^2(X,\Z)$-torsor by $Spin^c(X)$.

Recall that   the choice of a $Spin^c(4)$-structure $(\spm,\iota,\gamma)$
defines an
isomorphism between   the affine space ${\cal A}(\det \sp)$ of unitary
connections in
$\det\sp$ and the affine space of connections in $\Sigma^{\pm}$ which  lift
the
Levi-Civita connection in the bundle $\Lambda^2_{\pm}\simeq su(\spm)$. We
denote by
$\hat a\in{\cal A}(\Sigma)$ the connection corresponding to $a\in{\cal
A}(\det\sp)$.

The \underbar{Dirac} \underbar{operator} associated with the   connection
$a\in{\cal
A}(\det\sp)$ is the composition
$$\Dr_a:A^0(\spm)\textmap{\nabla_{\hat a}} A^1(\spm)\textmap{\gamma} A^0(\smp)
$$
of the covariant derivative $\nabla_{\hat a}$ in the bundles $\spm$ and the
Clifford
multiplication
$\gamma:\Lambda^1\otimes\spm\map \smp$.

Note that, in order to define the Dirac operator, one needs a Clifford map,
not only a
Riemannian metric and a pair of spinor bundles; this  will later become
important in
connection with transversality arguments.  The Dirac operator
$\Dr_a:A^0(\spm)\map
A^0(\smp)$ is an  elliptic first order operator with symbol
$\gamma:\Lambda^1\map\RSU(\spm,\smp)$. The  direct sum-operator
$\Dr_a:A^0(\Sigma)\map A^0(\Sigma)$ on the  total spinor bundle is
selfadjoint  and its
square has the same symbol as the rough Laplacian $\nabla_{\hat
a}^*\nabla_{\hat a}$ on
$A^0(\Sigma)$.

The corresponding \ub{Weitzenb\"ock} \ub{formula}   is
$$\Dr_a^2=\nabla_{\hat a}^*\nabla_{\hat a} +\frac{1}{2} \Gamma(F_a)+
\frac{s}{4}\id_\Sigma\ ,
$$
where $F_a\in i A^2$ is the curvature of the connection $a$, and $s$
denotes the scalar
curvature of $(X,g)$  [LM].

To write down the \sw equations, we need the following notations: For a
connection
$a\in{\cal A}(\det\sp)$ we let $F_a^{\pm}\in i A^2_{\pm}$ be the (anti) selfdual
components of its curvature. Given a spinor $\Psi\in A^0(\sp)$, we denote by
$(\Psi\bar\Psi)_0\in A^0(\End_0(\spm))$ the trace free part of the Hermitian
endomorphism $\Psi\otimes\bar\Psi$.
Now fix a  $Spin^c(4)$-structure  $\tau=(\spm,\iota,\gamma)$ for $(X,g)$
and a closed
2-form $\beta\in A^2$.  The $\beta$-\underbar{twisted} \underbar{monopole}
\underbar{equations} for a pair  $(a,\Psi)\in {\cal A}(\det\sp)\times
A^0(\sp)$ are
$$\left\{\begin{array}{ccc}
\Dr_a\Psi&=&0\\
\Gamma(F_a^++2\pi i\beta^+)&=&(\Psi\bar\Psi)_0 \ .
\end{array}
\right.\eqno{(SW^\tau_\beta)}
$$
These $\beta$-twisted  \sw equations should not be regarded as
perturbations of the equations
$(SW^\tau_0)$ since later the cohomology class of $\beta$ will be fixed.
The twisted
equations arise naturally in connection with non-abelian monopoles (see
section 2.2).
Using the Weitzenb\"ock formula one gets easily

\begin{lm} Let $\beta$ be a closed 2-form  and  $(a,\Psi)\in {\cal
A}(\det\sp) \times
A^0(\sp)$. Then we have the identity
$$
\begin{array}{l}\nr\Dr_a\Psi\nr^2+\frac{1}{4}\nr (F_a^++2\pi
i\beta^+)-(\Psi\bar\Psi)_0\nr^2=\\ \\
\nr\nabla_{A_a}\Psi\nr^2+\frac{1}{4}\nr F_a^++2\pi
i\beta^+\nr^2   +\frac{1}{8}\nr\Psi\nr_{L^4}^4+
\int\limits_X((\frac{s}{4}\id_{\sp}-\Gamma(\pi i\beta^+))\Psi,\Psi)  .
\end{array}$$
\end{lm}
\begin{co} {} [W] On manifolds $(X,g)$ with non-negative scalar curvature
$s$ the only solutions
of $(SW^\tau_0)$  are pairs $(a,0)$ with $F_a^+=0$.
\end{co}
\subsection{\sw invariants for 4-manifolds with $b_+>1$}
Let $(X,g)$  be a closed oriented Riemannian 4-manifold, and let $\cg\in
Spin(X)$ be an
equivalence class of $Spin^c$-structures of Chern class $c$, represented by
the triple
$\tau=(\spm,\iota,\gamma)$. The \underbar{configuration} \underbar{space}
for   \sw
theory is the product ${\cal A}(\det \sp)\times A^0(\sp)$ on which the
\underbar{gauge}
\underbar{group} ${\cal G}:={\cal C}^{\infty}(X,S^1)$ acts by
$$f\cdot(a,\Psi):=(a-2f^{-1}df\ ,\ f\Psi)\ .
$$

Let ${\cal B}(c):=\qmod{{\cal A}(\det \sp)\times A^0(\sp)}{{\cal G}}$ be
the orbit space; it
depends -- up to homotopy equivalence -- only on the Chern class $c$. Since
the gauge group
acts freely in all points $(a,\Psi)$ with $\Psi\ne 0$, the open subspace
$${\cal B}(c)^*:= \qmod{{\cal A}(\det \sp)\times
\left(A^0(\sp)\setminus\{0\}\right)}
{{\cal G}}$$ is a classifying space for ${\cal G}$. It has the weak
homotopy  type of  a
product $K(\Z,2) \times K(H^1(X,\Z),1)$ of Eilenberg-Mac Lane spaces and
there is a natural
isomorphism
$$\nu:\Z[u] \otimes\Lambda^*(\qmod{H_1(X,\Z)}{\rm Tors}) \map H^*({\cal
B}(c)^*,\Z)\ ,$$
where the generator $u$ is  of degree 2.  The ${\cal G}$-action on ${\cal
A}(\det\sp)\times
A^0(\sp)$ leaves the subset $[{\cal A}(\det\sp)\times
A^0(\sp)]^{SW^\tau_\beta}$ of solutions of $(SW^\tau_\beta)$ invariant; the
orbit space
$${\cal W}^\tau_{\beta}:=\qmod{[{\cal A}(\det\sp)\times
A^0(\sp)]^{SW^\tau_\beta}}{{\cal G}}
$$
is the \underbar{moduli} \underbar{space} of \ub{$\beta$}-\ub{twisted}
\ub{monopoles}. It
depends, up to canonical isomorphism, only on the metric $g$, on the closed
2-form $\beta$,
and on the class $\cg\in Spin^c(X)$ [OT6].

Let ${{\cal W}^\tau_\beta}^*\subset{\cal W}^\tau_\beta$ be the open subspace of
monopoles with non-vanishing spinor-component; it can be described  as the
zero-locus of a section in a vector-bundle over ${\cal B}(c)^*$. The total
space of this bundle
is
$$\left[{\cal A}(\det \sp)\times
\left(A^0(\sp)\setminus\{0\}\right)\right]\times_{\cal
G}\left[i A^2_+\oplus A^0(\Sigma^-)\right]\ ,$$
and the section is induced by the
${\cal G}$-equivariant map 
$$SW^\tau_\beta:{\cal A}(\det
\sp)\times \left(A^0(\sp)\setminus\{0\}\right)\map  iA^2_+\oplus
A^0(\Sigma^-)
$$
given by the  equations $(SW^\tau_\beta)$.

Completing the configuration space and the gauge group with respect to suitable
Sobolev norms, we can identify ${{\cal W}^\tau_\beta}^*$ with the zero set
of a real
analytic Fredholm section in the corresponding Hilbert   vector bundle on
the Sobolev
completion of ${\cal B}(c)^*$, hence we can endow this moduli space with
the structure of a
finite dimensional real analytic space.  As in the instanton case, one has
a Kuranishi
description for  local models of the moduli space around a given point
$[a,\Psi]\in{\cal
W}^\tau_\beta$ in terms of the first two cohomology groups of the elliptic
complex
$$0\map i A^0\textmap{D_p^0} iA^1\oplus A^0(\Sigma^+)\textmap{D_p^1}
iA^2_+\oplus
A^0(\Sigma^-)\map 0 \eqno{({\cal C}_p)}
$$
obtained by linearizing in $p=(a,\Psi)$ the action of the gauge group and
the equivariant map
$SW^\tau_\beta$. The differentials of this complex are
$$D_p^0(f)=(-2df,f\Psi) ,$$ $$
D_p^2(\alpha,\psi)=\left(d^+\alpha-\Gamma^{-1}[(\Psi\bar\psi)_0+
(\psi\bar\Psi)_0],\Dr_a(\psi)+\gamma(\alpha)(\Psi)\right) ,
$$
and its \ub{index} $w_c$ depends only on the Chern class $c$ of the
$Spin^c$-structure $\tau$  and on the  characteristic classes of the base
manifold $X$:
$$w_c=\frac{1}{4}(c^2-3\sigma(X)-2e(X)) \ .
$$

The moduli space ${\cal W}^\tau_{\beta}$ is compact.  This follows, as in
[KM],   from the
following consequence of the
 Weitzenb\"ock formula and the maximum principle.
\begin{pr} (A priori ${\cal C}^0$-bound of the spinor component) If
$(a,\Psi)$ is  a solution of
$(SW^\tau_\beta)$ then
$$\sup|\Psi|^2\leq \max\left(0, \sup\limits_X(-s+|4\pi\beta^+|)\right)\ .
$$
\end{pr}

Moreover,  let $\tilde {\cal W}^\tau$ be
the moduli space of triples $(a,\Psi,\beta)\in {\cal A}(\det\sp)\times
A^0(\sp)\times Z^2_{DR}$  solving the Seiberg-Witten equations above  now
regarded  as
equations for $(a,\Psi,\beta)$.    Two  such triples define the same point
in $\tilde {\cal
W}^\tau$ if they are congruent modulo the gauge group ${\cal G}$ acting
trivially on the third
component.  Using the proposition above and arguments of [KM], one can
easily see that the
natural projection
$\tilde {\cal W}^\tau\textmap{p} Z^2_{DR}$ is proper. Moreover, one has the
following
transversality results:
\begin{lm}   After suitable Sobolev completions the following holds:\\
1.  [KM] The open subspace
$[\tilde {\cal W}^\tau]^*\subset  \tilde {\cal W}^\tau$ of points with
non-vanishing spinor
component is smooth.  \\
2. [OT6] For every  de Rham cohomology class $b\in H^2_{DR}(X)$ the moduli
space  $[\tilde {\cal
W}^\tau_\beta]^*:=[\tilde {\cal W}^\tau]^*\cap p^{-1}(b)$ is also smooth.

\end{lm}

Now let $c\in H^2(X,\Z)$ be a \ub{characteristic} \ub{element}, i.e. an
integral lift of
$w_2(X)$. A pair $(g,b)\in{\cal M}et(X)\times H^2_{DR}(X)$ consisting of a
Riemannian metric
$g$ on $X$ and a de Rham cohomology class $b$ is called \ub{$c$}-\ub{good}
when the
$g$-harmonic representant of $c-b$ is not $g$-antiselfdual. This condition
guarantees that
${\cal W}^\tau_\beta={{\cal W}^\tau_\beta}^*$ for  every $Spin^c$-structure
$\tau$ of Chern
class $c$ and every 2-form $\beta$ in $b$. Indeed, if $(a,0)$ would solve
$(SW^\tau_\beta)$, then the $g$-antiselfdual 2-form
$\frac{i}{2\pi}F_a-\beta$ would be the
$g$-harmonic representant of $c-b$.

In particular, using the transversality results above, one gets the following
\begin{thry} { } [OT6] Let $c\in H^2(X,\Z)$ be a characteristic element and
suppose
$(g,b)\in{\cal M}et(X) \times H^2_{DR}(X)$ is $c$-good. Let $\tau$ be a
$Spin^c$-structure
of Chern class $c$ on $(X,g)$, and $\beta\in b$ a general representant of
the cohomology class
$b$. Then   the moduli space ${\cal W}^\tau_\beta={{\cal W}^\tau_\beta}^*$
is a closed manifold of dimension
$w_c=\frac{1}{4}(c^2-3\sigma(X)-2e(X))$.
\end{thry}

Let us fix a maximal subspace $H^2_+(X,\R)$ of $H^2(X,\R)$ on which the
intersection form
is positive definite. The dimension $b_+(X)$ of such a subspace is the number of positive
eigenvalues of the intersection form.   The moduli space  ${\cal
W}^\tau_\beta$ can be oriented by the choice of an orientation of the line $\det
H^1(X,\R)\otimes \det H^2_+(X,\R)^{\vee}$.

Let
$[{\cal W}^\tau_\beta]_\ooo\in H_{w_c}({\cal B}(c)^*,\Z)$ be the
\ub{fundamental} \ub{class} associated with the choice of an orientation
$\oo$ of the line
$\det H^1(X,\R)\otimes \det  H^2_{+}(X)^{\vee}$.

The \ub{Seiberg}-\ub{Witten} \ub{form} associated with the data
$(g,b,\cg,\oo)$ is the
element $SW^{(g,b)}_{X,\ooo}(\cg)\in \Lambda^*(H^1(X,\Z))$ defined by
$$SW^{(g,b)}_{X,\ooo}(\cg)(l_1\wedge\dots\wedge
l_r):=\left\langle\nu(l_1)\cup\dots
\nu(l_r)\cup u^{\frac{w_c-r}{2}},[{\cal W}^\tau_\beta]_\ooo\right\rangle
$$
for decomposable elements $l_1\wedge\dots\wedge l_r$ with $r\equiv w_2$
(mod 2). Here
$\tau$ is a $Spin^c$-structure on $(X,g)$ representing the class $\cg\in
Spin(X)$, and
$\beta$ is a general form in the class $b$.

One shows, using again transversality arguments, that the Seiberg-Witten form
$SW^{(g,b)}_{X,\ooo}(\cg)$ is well defined, independent of the choices of
$\tau$ and
$\beta$. Moreover, if any two $c$-good pairs $(g_0,b_0)$, $(g_1,b_1)$ can
be joined by a
smooth path of $c$-good pairs, then $SW^{(g,b)}_{X,\ooo}(\cg)$ is also
independent of
$(g,b)$ [OT6].

Note that the condition "$(g,b)$ is not $c$-good" is of codimension
$b_+(X)$ for a fixed class
$c$. This means that for manifolds with $b_+(X)>1$ we have a well defined map
$$SW_{X,\ooo}:Spin^c(X)\map  \Lambda^* H^1(X,\Z)
$$
which associates to a class of $Spin^c$-structures $\cg$ the form
$SW^{(g,b)}_{X,\ooo}(\cg)$ for any $b\in H^2_{DR}(X)$ such that $(g,b)$ is
$c$-good.  This
map, which is functorial  with respect to orientation preserving
diffeomorphisms, is the
\ub{Seiberg}-\ub{Witten} \ub{invariant}.

Using the identity in Lemma 1.1.1, one can easily prove
\begin{re} {}[W] Let $X$ be an oriented  closed 4-manifold with $b_+(X)>1$.
Then the set of
classes
$\cg\in Spin^c(X)$ with nontrivial Seiberg-Witten invariant is finite.
\end{re}

In the special case $b_+(X)>1$, $b_1(X)=0$, $SW_{X,\ooo}$ is simply a function
$$SW_{X,\ooo}:Spin^c(X)\map  \Z\ .
$$

The values $SW_{X,\ooo}(\cg)\in\Z$ are refinements of the numbers
$n_c^\ooo$ defined by
Witten [W]. More precisely:
$$ n_c^{\ooo}=\sum\limits_\cg SW_{X,\ooo}(\cg)\ ,
$$
the summation being over all classes of $Spin^c$-structures $\cg$ of Chern
class $c$. It is easy
to see that the indexing set is a torsor for the subgroup $Tors_2
H^2(X,\Z)$ of 2-torsion classes
in $H^2(X,\Z)$.

The structure of the Seiberg-Witten invariants for manifolds with
$b_+(X)=1$ is more
complicated and will be  described in the next section.
\newpage
\subsection{The case $b_+=1$ and the wall crossing formula}

Let $X$ be a closed  oriented differentiable 4-manifold with $b_+(X)=1$. In
this situation the
Seiberg-Witten forms depend on a \ub{chamber} \ub{structure}:
Recall first that in the case $b_+(X)=1$  there is a natural map ${\cal
M}et(X)\map  \P(H^2_{DR}(X))$ which sends a metric $g$ to the line
$\R[\omega_+]\subset
H^2_{DR}(X)$, where $\omega_+$ is any non-trivial $g$-selfdual harmonic
form. Let $\bf
H$ be the hyperbolic space
$${\bf H}:=\{h\in H^2_{DR}(X)|\ h^2=1\}\ .
$$

This space has two connected components, and the choice of one of them
orients the lines
$\H^2_{+,g}(X)$ for all metrics $g$. Furthermore, having fixed a component
${\bf H}_0$ of
${\bf H}$, every metric defines a unique $g$-selfdual form $\omega_g$ of
length 1 with
$[\omega_g]\in {\bf H}_0$.
\dfigure 100mm by 130mm (hyp scaled 500 offset 0mm:)

Let $c\in H^2(X,\Z)$ be characteristic. The \ub{wall} associated with $c$
is the hypersurface
$$c^{\bot}:=\{(h,b)\in {\bf H}\times H^2_{DR}(X)|\ (c-b)\cdot h=0\}\ ,
$$
and the connected components of  $[{\bf H}\times H^2_{DR}(X)]\setminus
c^{\bot}$ are called
\ub{chambers} of type $c$.

Notice that the walls are non-linear. Each characteristic element $c$
defines precisely four
chambers of type $c$, namely
$$C_{{\bf H}_0,\pm}:=\{(h,b)\in {\bf H}\times H^2_{DR}(X)|\ \pm(c-b)\cdot h<0\}\ ,\ {\bf
H}_0\in\pi_0({\bf H})\ ,
$$

and each of these four chambers contains elements of the form
$([\omega_g],b)$ with
$g\in{\cal M}et(X)$.

Let $\oo_1$ be an orientation of $H^1(X,\R)$. The choice of $\oo_1$
together with the choice of a
component ${\bf H}_0\in\pi_0({\bf H})$ defines an orientation
$\oo=(\oo_1,{\bf H}_0)$ of
$\det(H^1(X,\R))\otimes\det (H^2_+(X,\R)^{\vee})$. Set
$$SW_{X,(\ooo_1,{\bf H}_0)}^{\pm}(\cg):=SW_{X,\ooo}^{(g,b)}(\cg)\ ,
$$
where $(g,b)$ is a pair such that $([\omega_g],b)$ belongs to the chamber
$C_{{\bf H}_0,\pm}$.
The map
$$SW_{X,(\ooo_1,{\bf H}_0)}:Spin^c(X)\map \Lambda^* H^1(X,\Z)\times
\Lambda^* H^1(X,\Z)
$$
which associates to a class $\cg$ of $Spin^c$-structures on the oriented
manifold $X$ the pair of
forms  $(SW^{+}_{X,(\ooo_1,{\bf H}_0)}(\cg),SW^{-}_{X,(\ooo_1,{\bf
H}_0)}(\cg))$ is the
\ub{Seiberg}-\ub{Witten} \ub{invariant} of   $X$   with respect to the
orientation
data $(\oo_1,{\bf H}_0)$. This invariant is functorial with respect to
orientation-preserving
diffeomorphisms and behaves as follows with respect to changes of the
orientation data:
$$SW_{X,(-\ooo_1,{\bf H}_0)}(\cg)=-SW_{X,(\ooo_1,{\bf H}_0)}(\cg)\ ,\ \
SW^{\pm}_{X,(\ooo_1,-{\bf H}_0)}(\cg)=-SW^{\mp}_{X,(\ooo_1,{\bf H}_0)}(\cg) \ .
$$

More important, however, is the fact that the difference
$$SW^{+}_{X,(\ooo_1,{\bf H}_0)}(\cg)-SW^{-}_{X,(\ooo_1,{\bf H}_0)}(\cg)
$$
is a topological invariant of the pair $(X,c)$. To be precise, consider the
element $u_c\in
\Lambda^2\left(\qmod{H_1(X,\Z)}{Tors}\right)$ defined by
$$u_c(a\wedge b):=\frac{1}{2}\langle a\cup b\cup c,[X]\rangle
$$
for elements $a, b\in H^1(X,\Z)$. The following \ub{universal}
\ub{wall}-\ub{crossing}
formula generalizes results of
 [W], [KM] and [LL].
\begin{thry} {}[OT6] (Wall crossing formula) Let $l_{\ooo_1}\in
\Lambda^{b_1} H^1(X,\Z)$ be the
generator defined by the orientation $\oo_1$, and let $r\geq 0$ with
$r\equiv w_c$ (mod 2). For
every $\lambda\in \Lambda^r\left(\qmod{H_1(X,\Z)}{Tors}\right)$ we have
$$[SW^{+}_{X,(\ooo_1,{\bf H}_0)}(\cg)-SW^{-}_{X,(\ooo_1,{\bf
H}_0)}(\cg)](\lambda)=
\frac{ (-1)^{\frac{b_1-r}{2}} }  {\frac{b_1-r}{2}}\langle \lambda\wedge
u_c^{\frac{b_1-r}{2}},l_{\ooo_1}\rangle
$$
when $r\leq\min(b_1,w_c)$, and the difference vanishes otherwise.
\end{thry}
We want to illustrate these results with the simplest possible example, the
projective plane.\\
\\
{\bf Example:} Let $\P^2$ be the complex projective plane, oriented as
complex manifold, and
denote by $h$ the first Chern class of ${\cal O}_{\P^2}(1)$. Since $h^2=1$,
the hyperbolic space
${\bf H}$ consists of two points ${\bf H}=\{\pm h\}$. We choose the
component ${\bf
H}_0:=\{h\}$ to define orientations.

An element $c\in H^2(X,\Z)$ is characteristic iff $c\equiv h$ (mod 2). In
the picture below we
have drawn (as vertical intervals) the two chambers
$$C_{{\bf H}_0,\pm}=\{(h,b)\in {\bf H}_0\times H^2_{DR}(\P^2)|\
\pm(c-b)\cdot h<0\}
$$
of type $c$, for every $c\equiv h$ (mod 2).
\dfigure 187mm by 186mm (plane scaled 350 offset 0mm:)

The set $Spin^c(\P^2)$ can be identified with the set $(2\Z+1)h$ of
characteristic elements
under the map which sends a $Spin^c$-structure  $\cg$ to its Chern class
$c$. The corresponding
virtual dimension is $w_c=\frac{1}{4}(c^2-9)$.  Note that, for any metric
$g$, the pair $(g,0)$ is
$c$-good for all characteristic elements $c$. Also recall that the
Fubini-Study metric $g$ is a
metric of positive scalar curvature which can be normalized such that
$[\omega_g]=h$. We can
now completely determine the \sw invariant  $SW_{\P^2,{\bf H}_0}$ using
three simple
arguments:\\
i) for $c=\pm h$ we have $w_c<0$, hence $SW^{\pm}_{X,{\bf H}_o}(c)=0$, by
the transversality
results of section 1.2. \\
ii) Let $c$ be a characteristic element with $w_c\geq 0$. Since the
Fubini-Study metric $g$ has
positive scalar curvature and $(g,0)$ is $c$-good, we have ${\cal
W}^\tau_{0}={{\cal
W}^\tau_{0}}^*=\emptyset$ by  Corollary 1.1.2. But this
moduli space can be used to compute   $SW^{\pm}_{\P^2,{\bf H}_0}(c)$ for   characteristic
elements $c$ with $\pm c\cdot h<0$. Thus we find $SW^{\pm}_{\P^2,{\bf
H}_0}(c)=0$ when
$w_c\geq 0$ and
$\pm c\cdot h<0$.\\
iii) The remaining values, $SW^{\mp}_{\P^2,{\bf H}_0}(c)=0$ for classes
with $w_c\geq 0$ and
$\pm c\cdot h<0$, are determined by the wall-crossing formula. Altogether we get
$$SW^+_{\P^2,{\bf H}_0}(c)=\left\{\begin{array}{lll}
1&{\rm if}& c\cdot h\geq 3\\
0&{\rm if}&c\cdot h<3\ ,
\end{array}\right.\  \
SW^-_{\P^2,{\bf H}_0}(c)=\left\{\begin{array}{lll}
-1&{\rm if}& c\cdot h\leq-3\\
0&{\rm if}&c\cdot h>-3\ .
\end{array}\right.
 $$

\section{Non-abelian Seiberg-Witten theory}
\subsection{G-monopoles}

Let $V$ be a Hermitian vector space, and let $U(V)$ be its group of unitary
automorphisms.
For any closed subgroup $G\subset U(V)$ which contains the central
involution $-\id_V$, we
define a new Lie group by
$$Spin^G(n):=Spin(n)\times_{\Z_2} G\ .
$$
By construction one has the following exact sequences:
$$1 \map Spin \map Spin^G\textmap{\delta}\qmod{G}{\Z_2}\map  1
$$
$$1\map  G\map  Spin^G\textmap{\pi} SO\map  1
$$
$$1\map  \Z_2 \map Spin^G\textmap{(\pi,\delta)} SO\times\qmod{G}{\Z_2}\map
1 \ ,
$$
where    $Spin\ (Spin^G, SO)$  denotes one of the groups  $Spin(n)\
(Spin^G(n),SO(n))$.

Given a $Spin^G$-principal bundle $P^G$ over a topological space, we form
the following
associated bundles:
$$\delta(P^G):=P^G\times_\delta(\qmod{G}{\Z_2})\ ,\ \
\G(P^G):=P^G\times_{\rm Ad} G\ ,\ \
\gr(P^G):= P^G\times_{\rm ad} \g \ ,
$$
where $\g$ stands for the Lie algebra of $G$. The group ${\cal G}$ of
sections  of the bundle
$\G(P^G)$ can be identified with the group of automorphism of $P^G$ over
the associated
$SO$-bundle $P^G\times_\pi SO$.

 Consider now an oriented manifold $(X,g)$, and
let $P_g$ be the $SO$-bundle of oriented $g$-orthonormal coframes. A
\ub{$Spin^G$}-\ub{structure} in
$P_g$ is a principal bundle morphism  $\sigma:P^G\map  P_g$ of type $\pi$
[KN]. An
isomorphism of
$Spin^G$-structures $\sigma$, $\sigma'$ in $P_g$ is a bundle isomorphism
$f:P^G\map  P'^G$
with $\sigma'\circ f=\sigma$. One shows that the data of a
$Spin^G$-structure in $(X,g)$ is
equivalent to the data  of a linear, orientation-preserving isometry
$\gamma:\Lambda^1\map  P^G\times_\pi \R^n$, which we call the \ub{Clifford}
\ub{map}
of the $Spin^G$-structure [T2].

In dimension 4, the spinor group $Spin(4)$ splits as
$$Spin(4)=SU(2)_+\times SU(2)_-=Sp(1)_+\times Sp(1)_-\ .$$
Using the
projections
$$p_{\pm}:Spin(4)\map  SU(2)_{\pm}
$$
one defines the adjoint bundles
$$\ad_{\pm}(P^G):=P^G\times_{\ad_{\pm}} su(2)\ .$$
Coupling $p_{\pm}$ with the natural representation of $G$ in $V$, we obtain
representations $\lambda_{\pm}:Spin^G(4)\map  U(\H_{\pm}\otimes_\C V)$ and
associated
\ub{spinor} \ub{bundles}
$$\Sigma^{\pm}(P^G):= P^G\times_{\lambda_{\pm}}(\H_{\pm}\otimes_{\C} V)\ .
$$

The Clifford map $\gamma:\Lambda^1\map  P^G\times_\pi\R^4$ of the
$Spin^G$-structure
yields identifications
$$\Gamma:\Lambda^2_{\pm}\map  \ad_{\pm}(P^G)\ .
$$
An interesting special case occurs, when $V$ is a Hermitian vector space
over the
quaternions, and $G$ is a subgroup of $Sp(V)\subset U(V)$. Then one can
define \ub{real}
\ub{spinor} \ub{bundles}
$$\Sigma^{\pm}_\R(P^G):= P^G\times_{\rho_{\pm}} (\H_{\pm}\otimes_\H V)\ ,
$$
associated with the representations
$$\rho_{\pm}:Spin^G(4) \map SO(\H_{\pm}\otimes_{\H} V)\ .
$$
{\bf Examples:}  Let $(X,g)$ be a closed  oriented Riemannian 4-manifold
with coframe
bundle $P_g$.\\ \\
1) $G=S^1$: A $Spin^{S^1}$-structure is just a $Spin^c$-structure as
described in Chapter 1
[T2].
\\
\\
2) $G=Sp(1)$:  $Spin^{Sp(1)}$-structures have been introduced in [OT5],
where they were
called $Spin^h$-structures.  The map which associates to a
$Spin^{Sp(1)}$-structure
$\sigma:P^h\map  P_g$ the first Pontrjagin class $p_1(\delta(P^{h})$ of the
associated
$SO(3)$-bundle $\delta(P^h)$, induces a bijection between the set of
isomorphism classes
of $Spin^{Sp(1)}$-structures in $(X,g)$ and the set
$$\{p\in H^4(X,\Z)|\ p\equiv w_2(X)^2\ {\rm (mod\ 4)}\}\ .$$
There is a 1-1 correspondence between isomorphism classes of
$Spin^{Sp(1)}$-structures
in $(X,g)$ and equivalence classes of triples $(\tau:P^{S^1}\map
P_g,E,\iota)$   consisting of a $Spin^{S^1}$-structure
$\tau$, a unitary vector bundle $E$ of rank 2, and an unitary  isomorphism
$\iota:\det\Sigma^+_\tau\map \det E.$  The equivalence relation is
generated by tensorizing
with Hermitian line bundles [OT5], [T2]. The associated  bundles are -- in
terms of these
data -- given by
$$\delta(P^h)=\qmod{P_E}{S^1}\ ,\  \ \G(P^h)=SU(E)\  ,  \ \gr(P^h)=su(E), $$
$$
\Sigma^{\pm}(P^h)=(\Sigma^{\pm}_\tau)^{\vee}\otimes E\ , \  \
\Sigma^{\pm}_\R(P^h)=\RSU(\Sigma^{\pm}_\tau,E)\ ,
$$
where $P_E$ denotes the principal $U(2)$-frame bundle of $E$.
\\ \\
3) $G=U(2)$:  In this case $\qmod{G}{\Z_2}$ splits as
$\qmod{G}{\Z_2}=PU(2)\times S^1$,
and we write $\delta$ in the form $(\bar\delta,\det)$. The map which
associates to a
$Spin^{U(2)}$-structure $\sigma:P^u\map P_g$ the characteristic classes
$p_1(\bar\delta(P^u))$, $c_1(\det P^u )$ identifies the set of isomorphism
classes of
$Spin^{U(2)}$-structures in $(X,g)$ with the set
$$\{(p,c)\in H^4(X,\Z)\times H^2(X,\Z)|\ p\equiv (w_2(X)+\bar c)^2\ {\rm
(mod\ 4)}\}\ .
$$
There is a 1-1 correspondence between isomorphism classes of
$Spin^{U(2)}$-structures in
$(X,g)$ and equivalence classes of pairs $(\tau:P^{S^1}\map  P_g,E)$
consisting of a
$Spin^{S^1}$-structure $\tau$  and a unitary vector bundle of rank 2. Again
the equivalence
relation is given by tensorizing with Hermitian line bundles [T2]. If
$\sigma:P^u\map P_g$
corresponds to the pair $(\tau:P^{S^1}\map  P_g,E)$, the associated bundles
are now
$$\bar\delta(P^u)=\qmod{P_E}{S^1}\ ,\ \ \det P^u = \det
[\Sigma^+_\tau]^{\vee}\otimes\det
E\ ,
$$
$$
\G(P^u)=U(E)\ ,\ \ \gr(P^u)=u(E)\ ,\ \
\Sigma^{\pm}(P^u)=(\Sigma^{\pm}_\tau)^{\vee}\otimes E\ .
$$
We will later also need the subbundles
$\G_0(P^u):=P^u\times_{\Ad}SU(2)\simeq SU(E)$ and
$\gr_0:=P^u\times_{\ad}su(2)\simeq su(E)$. The group of sections
$\Gamma(X,\G_0)$ can be
identified with the group of automorphisms of $P^u$ over
$P_g\times(P^u\times_{\det}
S^1)$.\\ \\

Consider now again a general $Spin^G$-structure $\sigma:P^G\map  P_g$ in
the 4-manifold
$(X,g)$. The spinor bundle $\Sigma^{\pm}(P^G)$ has $\H_{\pm}\otimes_\C V$
as standard
fiber, so that the standard fiber $su(2)_{\pm}\otimes \g$ of the bundle
$\ad_{\pm}(P^G)\otimes\gr(P^G)$ can be viewed as real subspace of
$\End(\H_{\pm}\otimes_\C V)$. We define  a quadratic map
$$\mu_{0G}:\H_{\pm}\otimes_\C V\map  su(2)_{\pm}\otimes\g
$$
by sending   $\psi\in\H_{\pm}\otimes_\C V$ to the orthogonal projection
$pr_{su(2)_{\pm}\otimes\g}(\psi\otimes\bar\psi)$ of the Hermitian  endomorphism
$(\psi\otimes\bar\psi)\in\End(\H_{\pm}\otimes_\C V)$. One can show that
$-\mu_{0G}$ is
the total (hyperk\"ahler) moment map for the $G$-action on the space
$\H_{\pm}\otimes_\C V$ endowed with the natural hyperk\"ahler structure
given by left
multiplication with quaternionic units [T2].

These maps give rise to quadratic bundle maps
$$\mu_{0G}: \Sigma^{\pm}(P^G) \map  \ad_{\pm}(P^G)\otimes\gr(P^G) \ .
$$

In the case $G=U(2)$ one can project $\mu_{0U(2)}$ on
$\ad_{\pm}(P^G)\otimes\gr_0(P^G)$
and gets a map
$$\mu_{00}: \Sigma^{\pm}(P^u) \map  \ad_{\pm}(P^u)\otimes\gr_0(P^u) \ .
$$

Note that a fixed $Spin^G$-structure $\sigma:P^G\map  P_g$ defines a
bijection between
connections $A\in{\cal A} (\delta(P^G))$ in $\delta(P^G)$ and connections
$\hat A\in{\cal
A}(P^G)$ in the $Spin^G$-bundle $P^G$ which lift the Levi-Civita connection
in $P_g$ via
$\sigma$. This follows immediately from the third exact sequence above. Let
$$\Dr_A:A^0(\Sigma^{\pm}(P^G))\map  A^0(\Sigma^{\mp}(P^G))
$$
be the associated \ub{Dirac} \ub{operator}, defined by
$$\Dr_A:A^0(\Sigma^{\pm}(P^G))\textmap{\nabla_{\hat A}}
A^1(\Sigma^{\pm}(P^G))\textmap{\gamma} A^0(\Sigma^{\mp}(P^G))\ .
$$
Here $\gamma:\Lambda^1\otimes\Sigma^{\pm}(P^G) \map\Sigma^{\mp}(P^G)$ is
the Clifford
multiplication corresponding to the embeddings $\gamma:\Lambda^1\map
P^G\times_\pi\R^4\subset\Hom_\C(\Sigma^{\pm}(P^G),\Sigma^{\mp}(P^G))$.

\begin{dt}
Let $\sigma:P^G\map P_g$ be a $Spin^G$-structure in the Riemannian manifold
$(X,g)$. The
\ub{$G$}-\ub{monopole} \ub{equations} for a pair $(A,\Psi)\in{\cal
A}(\delta(P^G))\times
A^0(\Sigma^+(P^G))$ are
$$\left\{\begin{array}{ccc}
\Dr_A\Psi&=&0\\
\Gamma(F_A^+)&=&\mu_{0G}(\Psi)\ .
\end{array}
\right. \eqno{(SW^\sigma)}$$
\end{dt}
The solutions of these equations will be called \ub{$G$}-\ub{monopoles}.
The symmetry
group of the
$G$-monopole equations is the gauge group
${\cal G}:=\Gamma(X,\G(P^G))$. If the Lie algebra of $G$ has a non-trivial
center
$z(\g)$, then one has  a family of ${\cal G}$-equivariant "twisted"
$G$-monopole equations
$(SW^\sigma_\beta)$  parameterized by $iz(\g)$-valued   2-forms $\beta\in
A^2 (iz(\g))$:
$$\left\{\begin{array}{ccc}
\Dr_A\Psi&=&0\\
\Gamma((F_A+2\pi i \beta)^+)&=&\mu_{0G}(\Psi) \ .
\end{array}
\right. \eqno{(SW^\sigma_\beta)}$$

We denote by ${\cal M}^\sigma$, respectively  ${\cal M}^\sigma_\beta$ the
corresponding
moduli spaces of solutions modulo the gauge group ${\cal G}$.\\

Since in the   case $G=U(2)$ there exists the splitting
$$\qmod{U(2)}{\Z_2}=PU(2)\times S^1\ ,$$
 the data of a connection in
$\delta(P^u)=\bar\delta(P^u)\times_X\det P^u $ is equivalent to the data of
a pair of
connections $(A,a)\in {\cal A}(\bar\delta(P^u))\times {\cal A}(\det P^u )$.
This can be used to
introduce  new important equations, obtained by fixing the abelian
connection $a\in{\cal
A}(\det P^u )$ in the
$U(2)$-monopole equations, and
regarding it as a parameter. One gets in this way the  equations
$$\left\{\begin{array}{ccc}
\Dr_{A,a}\Psi&=&0\\
\Gamma(F_A^+)&=&\mu_{00}(\Psi)
\end{array}
\right. \eqno{(SW^\sigma_a})$$
for a pair $(A,\Psi)\in{\cal A}(\bar\delta(P^u))\times A^0(\Sigma^+(P^u))$,
which will \  be
called \  the \ \hbox{\ub{$PU(2)$}-\ub{monopole}} \ub{equations}.    These
equations should
be regarded as a  twisted version of the quaternionic monopole equations
introduced in
[OT5],  which coincide in our present framework with the
$SU(2)$-monopole equations. Indeed, a $Spin^{U(2)}$-structure
$\sigma:P^u\map  P_g$ with
trivialized determinant line bundle can be regarded as
$Spin^{SU(2)}$-structure, and   the
corresponding quaternionic monopole equations  are $(SW^\sigma_\theta)$,
where $\theta$
is  the trivial connection in $\det P^u$.

The $PU(2)$-monopole equations are only
invariant under the group
${\cal G}_0:=\Gamma(X,\G_0) $ of automorphisms of $P^u$ over
$P_g\times_X\det P^u $. We
denote by ${\cal M}^\sigma_a$ the \ub{moduli} \ub{space} of
$PU(2)$-monopoles modulo
this gauge  group.
Note that  ${\cal M}^\sigma_a$ comes with  a natural
\ub{$S^1$}-\ub{action} given by the formula
$\zeta\cdot[A,\Psi]:=[A,\zeta^{\frac{1}{2}}\Psi]$.\\
\\
\ub{Comparing with other formalisms:}\\ \\
1. For $G=S^1$, $V=\C$ one recovers the original abelian Seiberg-Witten
equations and the
twisted abelian \sw equations of [L1], [Bru], [OT6].\\
2. For $G=S^1$, $V=\C^{\oplus k}$ one gets the so called "multimonopole
equations" studied by
J. Bryan and R. Wentworth [BW].\\
3. In the case $G=U(2)$, $V=\C^2$ one obtains the $U(2)$-monopole equations
which  were
studied in [OT1]  (see also chapter 3).\\
4. In the case of a $Spin$-manifold $X$  and $G=SU(2)$ the corresponding
monopole equations have been  studied from a physical point of view  in [LM].
5. If $X$ is simply connected, the $S^1$- quotient $\qmod{{\cal
M}^\sigma_a}{S^1}$ of a
moduli space of $PU(2)$-monopoles can be identified with a moduli space  of
"non-abelian
monopoles"  as defined in [PT].   Note that in the general
 non-simply connected case, one has to use our formalism.

\begin{re} Let $G=Sp(n)\cdot S^1\subset U(\C^{2n})$ be the Lie group of
transformations of
$\H^{\oplus n}$   generated by left multiplication with quaternionic
matrices in $Sp(n)$ and
by right multiplication with complex numbers of modulus 1. Then
$\qmod{G}{\Z_2}$   splits as
$PSp(n)\times S^1$.  In the same way as in the $PU(2)$-case one defines the $PSp(n)$-monopole
equations $(SW^\sigma_a)$ associated with a $Spin^{Sp(n)\cdot
S^1}(4)$-structure
$\sigma:P^G\map  P_g$ in $(X,g)$ and an abelian connection $a$ in the associated
$S^1$-bundle.
\end{re}

The solutions of the (twisted) $G$- and
$PU(2)$-monopole equations are the absolute minima of certain gauge
invariant functionals on
the corresponding configuration spaces ${\cal A}(\delta(P^G))\times
A^0(\Sigma^+(P^G))$ and
${\cal A}(\bar \delta(P^u))\times A^0(\Sigma^+(P^G))$.  The investigation
of these
non-abelian \sw functionals and of their associated sets of critical points
is the subject of
forthcoming thesis by A. M. Teleman [Te].

 For simplicity we describe here only
the case of non-twisted $G$-monopoles. The \ub{Seiberg}-\ub{Witten}
\ub{functional}
$SW^\sigma:{\cal A}(\delta(P^G))\times A^0(\Sigma^+(P^G))\map \R$
associated to a
$Spin^G$-structure is defined by
$$SW^\sigma(A,\Psi):=\nr\nabla_{\hat A}\Psi\nr ^2+\frac{1}{4}\nr
F_A\nr^2+\frac{1}{2}\nr\mu_{0G}(\Psi)\nr^2+\frac{1}{4}\int\limits_X s|\Psi|^2\ .
$$

The \ub{Euler}-\ub{Lagrange}  \ub{equations} describing general critical
points are
$$\left\{\begin{array}{ccc}
d_A^* F_A+ J(A,\Psi)&=&0\ \\
\Delta_{\hat A}\Psi+\mu_{0G}(\Psi)(\Psi)+\frac{1}{4} s\Psi&=&0 \ ,
\end{array}\right.$$
where the current $J(A,\Psi)\in A^1(\gr(P^G))$ is given by
   $\sqrt{32}$ times the orthogonal
projection  of the
$\End(\Sigma^+(P^G))$-valued 1-form $\nabla_{\hat A}\Psi\otimes\bar\Psi\in
A^1(\End(\Sigma^+(P^G)))$ onto $A^1(\gr(P^G))$.

In the abelian case $G=S^1$, $V=\C$, a closely related functional and the
corresponding
Euler-Lagrange equations have been investigated in [JPW].

\subsection{Moduli spaces of $PU(2)$-monopoles}

We retain the notations of the previous section. Let $\sigma:P^u\map  P_g$ be a
$Spin^{U(2)}$-structure in a closed oriented Riemannian 4-manifold $(X,g)$,
and let $a\in{\cal
A}(\det P^u )$ be a fixed  connection. The $PU(2)$-monopole equations
$$\left\{\begin{array}{ccc}
\Dr_{A,a}\Psi&=&0\\
\Gamma(F_A^+)&=&\mu_{00}(\Psi)
\end{array}
\right. \eqno{(SW^\sigma_a})$$
associated with these data are invariant under the action of the gauge
group ${\cal G}_0$, and
hence give  rise to a closed subspace ${\cal M}^\sigma_a\subset {\cal
B}(P^u)$ of the orbit
space ${\cal B}(P^u):=\qmod{{\cal A}(\bar\delta(P^u))\times
A^0(\Sigma^+(P^u))}{{\cal G}_0}$.

The moduli space ${\cal M}^\sigma_a$ can be endowed with the structure of a
ringed space
with \ub{local} \ub{models} constructed by the well-known Kuranishi method
[OT1],
[OT5], [DK], [LT]. More precisely: The linearization of the
$PU(2)$-monopole equations in a
solution
$p=(A,\Psi)$ defines an  \ub{elliptic}  \ub{deformation}  \ub{complex}
$$0\rightarrow A^0(\gr_0(P^u))\stackrel{D_p^0}{\rightarrow}
A^1(\gr_0(P^u))\oplus
A^0(\Sigma^+(P^u))\stackrel{D^1_p}{\rightarrow}A^2_+(\gr_0(P^u))\oplus
A^0(\Sigma^-(P^u))\rightarrow 0
$$
whose differentials are given by $D_p^0(f)=(-d_A f,f\Psi)$ and
$$D^1_p(\alpha,\psi)=(d^+_A\alpha-\Gamma^{-1}[m(\psi,\Psi)+m(\Psi,\psi)],\Dr
_{A,a}\psi+
\gamma(\alpha) \Psi)\ .$$
Here $m$ denotes the sesquilinear map associated with the quadratic map
$\mu_{00}$.  Let
$\H^i_p$, $i=0,1,2$ denote the harmonic spaces of the elliptic   complex
above. The
stabilizer ${\cal G}_{0p}$ of the point $p\in {\cal A}(\bar\delta(P^u))\times
A^0(\Sigma^+(P^u))$ is a finite dimensional Lie group, isomorphic to a
closed subgroup of
$SU(2)$, which acts in a natural way on the spaces $\H^i_p$.
\begin{pr}{ }[OT5], [T2] For every point $p\in {\cal M}^\sigma_a$ there exists a
neighborhood
$V_p\subset {\cal M}^\sigma_a$, a  ${\cal G}_{0p}$-invariant neighborhood
$U_p$ of
$0\in\H^1_p$, an
${\cal G}_{0p}$-equivariant map $K_p:U_p\map  \H^2_p$  with $K_p(0)=0$,
$dK_p(0)=0$,
and an isomorphism of ringed spaces
$$V_p\simeq\qmod{Z(K_p)}{{\cal G}_{0p}}
$$
sending $p$ to $[0]$.
 \end{pr}
The local isomorphisms $V_p\simeq\qmod{Z(K_p)}{{\cal G}_{0p}}$ define  the
structure of a
smooth manifold on the open subset
$${\cal M}^\sigma_{a,\rm reg}:=\{[A,\Psi]\in {\cal M}^\sigma_a|\ {\cal G}_{0p}=\{1\},\
\H^2_p=\{0\}\}\ ,
$$
and a real analytic orbifold structure in the open set of points $p\in{\cal
M}^\sigma_a$
with
${\cal G}_{0p}$ finite. The dimension of ${\cal M}^\sigma_{a,\rm reg}$
coincides with the
\ub{expected}
\ub{dimension} of the $PU(2)$-monopole moduli space, which is given by
the index
$\chi(SW^\sigma_a)$ of the elliptic deformation complex:
$$
\chi(SW^\sigma_a)=\frac{1}{2}(-3p_1(\bar\delta(P^u))+c_1(\det P^u
)^2)-\frac{1}{2}
(3e(X)+4\sigma(X))\ .
$$

Our next  goal is to describe the fixed point set of the $S^1$-action on
${\cal M}^\sigma_a$
introduced above.

First consider the closed subspace ${\cal D}(\bar\delta(P^u))\subset{\cal
M}^\sigma_a$ of
points of the form
$[A,0]$. It can be identified with the Donaldson moduli space of
anti-selfdual connections
in the $PU(2)$-bundle $\bar\delta(P^u)$  modulo the gauge group
${\cal G}_0$. Note however, that if $H^1(X,\Z_2)\ne \{0\}$,
${\cal D}(\bar\delta(P^u))$ does not coincide   with the usual moduli space of
$PU(2)$-instantons in $\bar\delta(P^u)$  but is  a finite cover of it.

The stabilizer ${\cal G}_{0p}$ of a Donaldson  point $(A,0)$ contains
always $\{\pm\id\}$,
hence  ${\cal M}^\sigma_a$ has at least $\Z_2$-orbifold singularities in
the points of ${\cal
D}(\bar\delta(P^u))$.

Second consider   $S^1$ as a subgroup of $PU(2)$ via the standard
embedding $S^1\ni\zeta\longmapsto [\left(\matrix{\zeta&0\cr 0&1}\right)]\in
PU(2)$.  Note
that any $S^1$-reduction $\rho:P\map \bar\delta(P^u)$  of  $\bar\delta(P^u)$
defines a reduction $\tau_\rho:P^\rho:=P^u\times_{\bar\delta(P^u)}P\map
P^u\stackrel{\sigma}{\map} P_g$ of the
$Spin^{U(2)}$-structure  \ $\sigma$ to a $Spin^{S^1\times S^1}$-structure,
hence a
pair of \ $Spin^c$-structures \  $\tau^i_\rho:P^{\rho_i}\map P_g$. One has
natural
isomorphisms
$$\det P^{\rho_1} \otimes \det P^{\rho_2} =(\det P^u)^{\otimes 2}\ , \ \
\det P^{\rho_1} \otimes (\det P^{\rho_2})^{-1}=P^{\otimes 2}\ ,$$
 and natural embeddings
$\Sigma^{\pm}(P^{\rho_i})\map\Sigma^{\pm}(P^u)$ induced by the   bundle
morphism $P^{\rho_i}\map P^u$.   A pair $(A,\Psi)$ will be called
\ub{abelian} if  it lies in the
image of ${\cal A}(P)\times A^0(\Sigma^+(P^{\rho_1}))$ for a suitable
$S^1$-reduction $\rho$
of $\bar\delta(P^u))$.

\begin{pr} The fixed point set of the $S^1$-action on ${\cal M}^\sigma_a$
is the union of the
Donaldson locus ${\cal D}(\bar\delta(P^u))$ and the locus of abelian
solutions. The latter can be
identified with the  disjoint union $\coprod\limits_\rho{\cal
W}^{\tau_\rho^1}_{\frac{i}{2\pi}
F_a}$, where the  union is over all $S^1$-reductions of  the $PU(2)$-bundle
$\bar\delta(P^u)$.
\end{pr}

This result suggests to use the $S^1$-quotient of   ${\cal
M}^\sigma_a\setminus({\cal M}^\sigma_a)^{S^1}$  for the comparison of Donaldson
invariants and (twisted) Seiberg-Witten invariants, as explained in [OT5].

Note that only using
moduli spaces  ${\cal M}^\sigma_\theta$ of quaternionic monopoles one gets,
by the
proposition above,  moduli spaces of \ub{non}-\ub{twisted}   abelian
monopoles  in the
fixed point locus of the $S^1$-action.  This was one of the motivations for
studying the
quaternionic monopole equations in [OT5]. There it has been shown that one
can use the
moduli spaces of quaternionic monopoles to relate  certain $Spin^c$-polynomials
to the original non-twisted \sw invariants.

The remainder    of this section is devoted to the description of the
\ub{Uhlenbeck}
\ub{compactification}   of the moduli spaces of $PU(2)$-monopoles [T3].

First of all, the Weitzenb\"ock formula and the maximum principle
yield  a bound on the spinor component, as in the abelian case.   More
precisely, one has the
a priori estimate
$$\sup\limits_X |\Psi|^2\leq
C_{g,a}:=\max\left(0,C\sup(-\frac{s}{2}+|F_a^+|)\right)
$$
on the space of solutions of  $(SW^\sigma_a)$, where $C$ is a universal
positive constant.

The construction of the Uhlenbeck compactification of ${\cal M}^\sigma_a$
is based, as in the
instanton case, on the following three essential results.\\
1.  A \ub{compactness} theorem for the subspace of solutions with suitable
bounds on the
 curvature of the connection component.\\
2. A \ub{removable} \ub{singularities} theorem.\\
3.   {Controlling} \ub{bubbling} phenomena for  an arbitrary sequence of
points in  the
moduli space
${\cal M}^\sigma_{a}$.  \\ \\ \\
1.  \ub{A compactness result}.

\begin{thry}  There exists a positive number $\delta>0$ such that for every
oriented
Riemannian manifold $(\Omega,g)$ endowed with a
$Spin^{U(2)}(4)$-structure $\sigma:P^u\map  P_g$  and a fixed    connection
$a\in{\cal A} (\det P^u )$, the following holds:

If $(A_n,\Psi_n)$ is a sequence of solutions of
$(SW^\sigma_{a})$, such that any point
$x\in\Omega$ has a geodesic ball neighborhood
$D_x$ with
$$\int_{D_x}| F_{A_n}|^2< \delta^2 $$
for all large enough $n$, then there is a subsequence
$(n_m)\subset\N$ and gauge transformations $f_{m}\in{\cal G}_0$ such that
$f_m^*\left(A_{m_n},\Psi_{m_n}\right)$ converges in the ${\cal
C}^{\infty}$-topology on
$\Omega$.
\end{thry}

 2.  \ub{Removable singularities}\\

Let $g$ be a metric on the 4-ball  $B$, and let
$$\sigma:P^u=B\times Spin^{U(2)}(4)\map P_g\simeq B\times SO(4)$$
 be a $Spin^{U(2)}$-structure in $(B,g)$. Fix $a\in iA^1_B$  and put
$B^{\bullet}:=B\setminus\{0\}$,
$\sigma^{\bullet}:=\sigma|_{B^{\bullet}}$.
\begin{thry} Let $(A_0,\Psi_0)$ be a solution of  the equations
$(SW^{\sigma^{\bullet}}_{a})$  on the punctured ball such that
$$\nr F_{A_0}\nr_{L^2}^2<\infty \ .
$$
Then there exists a solution $(A,\Psi)$ of $(SW^\sigma_{a})$ on $B$ and a gauge
transformation $f\in{\cal C}^{\infty}(B^{\bullet},SU(2))$ such that
$f^*(A|_{B^{\bullet}},\Psi|_{B^{\bullet}})=(A_0,\Psi_0)$.
\end{thry}

 3. \ub{Controlling bubbling phenomena}\\

The main point is that the selfdual components $F_{A_n}^+$ of the
curvatures of a  sequence
of solutions $([A_n,\Psi_n])_{n\in\N}$ in ${\cal M}^\sigma_{a}$ cannot bubble.
\begin{dt} Let   $\sigma: P^u\map P_g$ be a $Spin^{U(2)}$-structure in
$(X,g)$ and fix
$a\in{\cal A}(\det P^u)$. An \ub{ideal}
\ub{monopole} of type $(\sigma,a)$ is a pair
$([A',\Psi'],\{x_1,\dots,x_l\})$ consisting of a
point $[A',\Psi']\in{\cal M}^{\sigma'_l}_{a}$, where  $\sigma'_l:P'^u\map
P_g$ is a
$Spin^{U(2)}$-structure  satisfying
$$\det P'^u =\det P^u \ ,\ \  p_1(\bar\delta (P'^u)
)=p_1(\bar\delta(P^u))+4l\ ,$$
 and $\{x_1,\dots,x_l\}\in S^l(X)$. The set of ideal monopoles of type
$(\sigma,a)$ is
$$I{\cal M}^\sigma_{a}:=\coprod\limits_{l\geq 0} {\cal
M}^{\sigma'_l}_{a}\times S^l  X\ .
$$
\end{dt}
\begin{thry} There exists a metric topology on $I{\cal M}^\sigma_{a}$ such
that the moduli space ${\cal M}^\sigma_{a}$ becomes an open subspace
with compact closure $\overline{{\cal M}^\sigma_{a}}$.
\end{thry}

\pf (sketch) Given a sequence $([A_n,\Psi_n])_{n\in\N}$ of points in ${\cal
M}^\sigma_{a}$, one finds a subsequence
$([A_{n_m},\Psi_{n_m}])_{m\in\N}$, a finite set of points $S\subset X$,
and gauge transformations $f_m$ such that
$(B_m,\Phi_m):=f_m^*(A_{n_m},\Psi_{n_m})$ converges on $X\setminus S$ in the
${\cal C}^{\infty}$-topology to a solution
$(A_0,\Psi_0)$. This follows from the compactness theorem above, using the
fact that
the total volume of the sequence of measures $|F_{A_n}|^2$ is bounded. The
set $S$
consists of points in which the measure $|F_{A_{n_m}}|^2$ becomes
concentrated as $m$
tends to infinity.

By the Removable Singularities theorem, the solution $(A_0,\Psi_0)$
extends after gauge
transformation to a solution $(A,\Psi)$ of $(SW^{\sigma'}_{a})$ on $X$, for
a possibly
different $Spin^{U(2)}$-structure $\sigma'$ with the same determinant line
bundle. The
curvature of $A$ satisfies
$$|F_A|^2=\lim\limits_{n\rightarrow \infty}|F_{A_{n_m}}|^2-8\pi^2\sum_{x\in S}
\lambda_x\delta_x
$$
where $\delta_x$ is the Dirac measure of the point $x$. Now it remains to
show that
the $\lambda_x$'s are natural numbers and
$\sum\limits_{x\in
S}\lambda_x=\frac{1}{4}(p_1(\bar\delta(P'^u))-p_1(\bar\delta(P^u)))$. This
follows as in
the instanton case, if one uses the fact that the measures $|F_{A_{n_m}}^+|^2$ cannot
bubble in the points   $x\in S$ as $m\rightarrow\infty$ and that the integral of
$|F_{A_n}^-|^2-|F_{A_n}^+|^2$ is a topological invariant of
$\bar\delta(P^u)$. In this way one gets an
ideal monopole $\mg:=([A,\Psi],\{\lambda_1x_1,\dots,\lambda_k x_k\})$ of type
$(\sigma,a)$. With respect to a suitable   topology on the space of ideal
monopoles,
one has
$\lim\limits_{m\rightarrow\infty}[A_{n_m},\Psi_{n_m}]=\mg$.
\qed

\section{Seiberg-Witten theory and K\"ahler geometry}
\subsection{Monopoles on K\"ahler surfaces}

Let $(X,J,g)$ be an almost Hermitian surface with associated K\"ahler form
$\omega_g$.
We denote by $\Lambda^{pq}$ the bundle of $(p,q)$-forms on $X$ and by
$A^{pq}$ its space
of sections. The Hermitian structure defines an orthogonal decomposition
$$\Lambda^2_+\otimes\C=
\Lambda^{20}\oplus\Lambda^{02}\oplus \Lambda^{00}\omega_g
$$
and a canonical $Spin^c$-structure $\tau$. The spinor bundles of $\tau$ are
$$\Sigma^+=\Lambda^{00}\oplus \Lambda^{02}\ ,\ \ \Sigma^-=\Lambda^{01}\ ,
$$
and the Chern class of $\tau$ is the first Chern class
$c_1(T^{10}_J)=c_1(K_X^{\vee})$ of the
complex tangent bundle. The complexification of the canonical Clifford map
$\gamma$ is
the standard isomorphism
$$\gamma:\Lambda^1\otimes\C \map
\Hom(\Lambda^{00}\oplus\Lambda^{02},\Lambda^{01})\ ,\ \
\gamma(u)(\varphi+\alpha)=\sqrt{2}(\varphi u^{01}-i\Lambda_g
u^{10}\wedge\alpha)\ ,
$$
and the induced isomorphism $\Gamma:\Lambda^{20}\oplus \Lambda^{02}\oplus
\Lambda^{00}\omega_g\map\End_0(\Lambda^{00}\oplus \Lambda^{02})$ acts by
$$(\lambda^{20},\lambda^{02},f\omega_g)\stackrel{\Gamma}{\longmapsto}
2\left[\matrix{-if&-*(\lambda^{20}\wedge\cdot)\cr
\lambda^{02}\wedge\cdot&if\cr}\right]\in\End_0(\Lambda^{00}\oplus\Lambda^{02
}).$$

Recall from section 1.1 that the set $Spin^c(X)$  of equivalence classes of
$Spin^c$-structures in $(X,g)$ is a $H^2(X,\Z)$-torsor. Using the class of
the canonical
$Spin^c$-structure $\cg:=[\tau]$ as base point, $Spin^c(X)$ can be
identified with the set of
isomorphism classes of $S^1$-bundles: When $M$ is an $S^1$-bundle with
$c_1(M)=m$, the
$Spin^c$-structure $\tau_m$ has spinor bundles $\Sigma^{\pm}\otimes M$ and
Chern class
$2c_1(M)-c_1(K_X)$.  Let $\cg_m$ be the class of $\tau_m$.

Suppose now that $(X,J,g)$ is K\"ahler, and let $k\in{\cal A}(K_X)$ be the
Chern connection
in the canonical line bundle. In order to write the (abelian) \sw equations
associated with
the $Spin^c$-structure $\tau_m$ in a convenient form, we make the variable
substitution
$a=k\otimes e^{\otimes 2}$ for a connection $e\in{\cal A}(M)$ in the
$S^1$-bundle $M$,
and we write the spinor $\Psi$ as a sum $\Psi=\varphi+\alpha\in
A^0(M)\oplus A^{02}(M)$.
\begin{lm} {} [W], [OT6] Let $(X,g)$ be a K\"ahler surface,  $\beta\in
A^{11}_\R$ a closed
real (1,1)-form in the de Rham cohomology class $b$, and let $M$ be a
$S^1$-bundle with
$(2c_1(M)-c_1(K_X)-b)\cdot[\omega_g]<0$. The pair $(k\otimes e^{\otimes
2},\varphi+\alpha)\in{\cal A}(\det(\Sigma^+\otimes M))\times
A^0(\Sigma^+\otimes M)$
solves the equations $(SW^{\tau_m}_{X,\beta})$ iff $\alpha=0$,
$F_e^{20}=F_e^{02}=0$,
$\bar\partial_e\varphi=0$, and
$$i\Lambda_g
F_e+\frac{1}{4}\varphi\bar\varphi+(\frac{s}{2}-\pi\Lambda_g\beta)=0\ .
\eqno{(*)}$$
\end{lm}

Note that the conditions $F_e^{20}=F_e^{02}=0$, $\bar\partial_e\varphi=0$
mean that $e$
is the Chern connection of a holomorphic structure in the Hermitian line
bundle $M$ and
that $\varphi$ is a holomorphic section with respect to this holomorphic
structure. Integrating
the relation $(*)$ and using the inequality in the hypothesis,  one    sees
that $\varphi$
cannot vanish.

To interpret the  condition
$(*)$ consider an arbitrary real valued function function
$t:X\map\R$, and let
$$m_t:{\cal A}(M)\times A^0(M) \map iA^0
$$
be the map defined by
$$m_t(e,\varphi):=\Lambda_gF_e-\frac{i}{4}\varphi\bar\varphi +it \ .
$$

It easy to see that (after suitable Sobolev completions) ${\cal A}(M)\times
A^0(M)$ has a
natural symplectic structure, and that $m_t$ is a \ub{moment} \ub{map} for the action of the
gauge group ${\cal G}={\cal C}^{\infty}(X,S^1)$. Let ${\cal G}^\C={\cal
C}^{\infty}(X,\C^*)$
be the complexification of ${\cal G}$, and let ${\cal H}\subset A(M) \times
A^0(M)$ be the
closed set
$${\cal H}:=\{(e,\varphi)\in A(M) \times A^0(M)|\
F_e^{02}=0,\ \bar\partial_e \varphi =0\}
$$
of integrable pairs. For any function $t$ put
$${\cal H}_t:=\{(e,\varphi)\in {\cal H}|\  {\cal
G}^\C(e,\varphi)\cap m_t^{-1}(0) \ne\emptyset\}\ .
$$
Using a general principle in the theory of symplectic quotients, which also
holds in our
infinite dimensional framework,  one can prove that the
${\cal G}^\C$-orbit of a point
$(e,\varphi)\in{\cal H}_t$ intersects the zero set $m_t^{-1}(0)$ of the
moment map
$m_t$ precisely along a
${\cal G}$-orbit.
\dfigure 149mm by 145mm (moment scaled 600 offset 0mm:)
In other words, there is a natural bijection of quotients
$$\qmod{[m_t^{-1}(0)\cap{\cal H}]}{{\cal G}}\simeq \qmod{{\cal H}_t}{{\cal
G}^\C}\ .
\eqno{(1)}
$$

Now take $t:=- (\frac{s}{2} -\pi\Lambda_g\beta)$ and suppose again that the
assumptions
in the proposition hold. We have seen that
$m_t^{-1}(0)\cap{\cal H}$ cannot contain pairs of the form $(e,0)$, hence
${\cal G}$
(${\cal G}^\C$) acts freely in $m_t^{-1}(0)\cap{\cal H}$ (${\cal H}_t$).
Using this fact one can show that ${\cal H}_t$ is open in the space ${\cal
H}$ of
integrable pairs, and endowing the two quotients in (1) with the natural
real analytic
structures, one   proves that (1) is a real analytic isomorphism. By the
proposition above, the first quotient is precisely the moduli space  ${\cal
W}^{\tau_m}_{\beta}$. The second quotient is a complex-geometric object,
namely an
open subspace in the moduli space of simple holomorphic pairs $\qmod{{\cal
H}\cap\{\varphi\ne 0\}}{{\cal G}^\C}$. A point in this moduli space can be
regarded as an
isomorphism class of pairs $({\cal M},\varphi)$ consisting of a holomorphic
line bundle
${\cal M}$ of topological type $M$, and a holomorphic section in ${\cal M}$.
Such a pair defines a point in
$\qmod{{\cal H}_t}{{\cal G}^\C}$ if and only if ${\cal M}$ admits a
Hermitian metric $h$
satisfying the equation
$$i\Lambda F_h+\frac{1}{4}\varphi\bar\varphi^h=t \ .\eqno{(V_t)}
$$
This  equation -- for the unknown metric  $h$ -- is the \ub{vortex}
\ub{equation} associated
with the function $t$: it is solvable iff the \ub{stability} \ub{condition}
$$(2c_1(M)-c_1(K_X)-b)\cdot [\omega_g]<0
$$
is fulfilled. Let ${\cal D}ou(m)$ be the Douady space of effective divisors
$D\subset X$
with $c_1({\cal O}_X(D))=m$. The map $Z:\qmod{{\cal H}_t}{{\cal G}^\C}\map {\cal
D}ou(m)$ which associates to an orbit $[e,\varphi]$ the zero-locus
$Z(\varphi)\subset X$ of
the holomorphic section $\varphi$ is an isomorphism of complex spaces.

Putting everything together, we have the following interpretation for the
monopole
moduli spaces ${\cal W}^{\tau_m}_\beta$ on K\"ahler surfaces.
\begin{thry} {} [OT1], [OT6] Let $(X,g)$ be a compact  K\"ahler surface,
and let $\tau_m$ be
the
$Spin^c$-structure defined by the $S^1$-bundle $M$. Let $\beta\in
A^{11}_\R$ be a closed
2-form representing the de Rham cohomology class $b$ such that
$$(2c_1(M)-c_1(K_X)-b)\cdot [\omega_g]<0\ (>0)\ .
$$
If $c_1(M)\not\in NS(X)$, then ${\cal W}^{\tau_m}_\beta=\emptyset$. When
$c_1(M)\in NS(X)$, then there is a natural real analytic isomorphism
$${\cal W}^{\tau_m}_\beta\simeq {\cal D}ou(m)\ \ ({\cal D}ou(c_1(K_X)-m)) \ .
$$
\end{thry}

A moduli space ${\cal W}^{\tau_m}_\beta\ne\emptyset$ is smooth at the point
corresponding to $D\in{\cal D}ou(m)$ iff $h^0({\cal O}_D(D))=\dim_D{\cal
D}ou(m)$. This
condition is always satisfied when $b_1(X)=0$. If ${\cal W}^{\tau_m}_\beta$
is smooth at a point corresponding to $D\in{\cal D}ou(m)$, then it has the
expected
dimension in this point iff $h^1({\cal O}_D(D))=0$.

The natural isomorphisms  ${\cal W}^{\tau_m}_\beta\simeq {\cal D}ou(m)$
respects the
orientations induced by the complex structure  of $X$ when
$(2c_1(M)-c_1(K_X)-b)\cdot
[\omega_g]<0$. If  $(2c_1(M)-c_1(K_X)-b)\cdot [\omega_g]>0$, then the isomorphism
${\cal W}^{\tau_m}_\beta\simeq {\cal D}ou(c_1(K_X)-m)$ multiplies the complex
orientations by $(-1)^{\chi(M)}$ [OT6].\\ \\
{\bf Example:} Consider again the complex  projective plane $\P^2$, polarized by
$h=c_1({\cal O}_{\P^2}(1))$.  The expected dimension of  ${\cal
W}^{\tau_m}_\beta$ is
$m(m+3h)$. The theorem above yields the following explicit description of the
corresponding moduli spaces:
$${\cal W}_{\P^2,\beta}^{\tau_m}\simeq\left\{\begin{array}{cll}
|{\cal O}_{\P^2}(m)|&{\rm if}&(2m+3h-[\beta])\cdot h<0\\
|{\cal O}_{\P^2}(-(m+3))|&{\rm if}&(2m+3h-[\beta])\cdot h>0\ .
\end{array}\right.
$$

E. Witten has shown [W] that on K\"ahlerian surfaces $X$ with geometric
genus $p_g>0$
all non-trivial \sw invariants $SW_{X,\ooo}(\cg)$ satisfy $w_c=0$.

In the case of K\"ahlerian surfaces with $p_g=0$ one has a different
situation. Suppose for
instance  that $b_1(X)=0$. Choose    the standard orientation $\oo_1$ of
$H^1(X,\R)=0$ and    the component ${\bf H}_0$ containing  K\"ahler classes
to orient the
moduli spaces of monopoles.  Then, using the previous theorem and the
wall-crossing
formula, we get:
\begin{pr} Let $X$ be a K\"ahler surface with $p_g=0$ and $b_1=0$. If $m\in
H^2(X,\Z)$
satisfies  $m(m-c_1(K_X))\geq 0$, i.e. the expected dimension
$w_{2m-c_1(K_X)}$ is
non-negative, then
$$SW^+_{X,{\bf H}_0}(\cg_M)=\left\{\begin{array}{ccc} 1&{\rm if} &{\cal
D}ou(m)\ne\emptyset\\
0 &{\rm if} &{\cal D}ou(m)=\emptyset \ ,
\end{array}\right.
$$
$$SW^-_{X,{\bf H}_0}(\cg_M)= \left\{\begin{array}{ccc} 0&{\rm if} &{\cal
D}ou(m)\ne\emptyset\\
-1 &{\rm if} &{\cal D}ou(m)=\emptyset \ .
\end{array}\right.
$$
\end{pr}
\vspace{7mm}

Our next  goal is to show that the $PU(2)$-monopole equations on a K\"ahler
surface can
be analyzed in a similar way. This analysis yields a complex geometric
description  of the
moduli spaces whose $S^1$-quotients   give formulas relating the
Donaldson invariants to the \sw invariants. If the base is projective, one
also has an
algebro-geometric interpretation [OST], which leads to explicitly
computable examples of
moduli spaces of $PU(2)$-monopoles [T3]. Such examples are important,
because they illustrate
the general mechanism for proving the relation between the two theories,
and help to
understand the geometry of the ends of the moduli spaces in the more
difficult ${\cal
C}^{\infty}$-category.\\

Recall that, since $(X,g)$ comes with a canonical $Spin^c$-structure
$\tau$, the data of
of a $Spin^{U(2)}$-structure in $(X,g)$ is equivalent to the data of a
Hermitian bundle $E$
of rank 2. The bundles of the corresponding $Spin^{U(2)}$-structure
$\sigma: P^u\map
P_g$ are given by $\bar\delta(P^u)=\qmod{P_E}{S^1}$, $\det P^u =\det E
\otimes K_X$,
and $\Sigma^{\pm}(P^u)=\Sigma^{\pm}\otimes E\otimes K_X$.

Suppose that $\det P^u $ admits an \ub{integrable} connection $a\in{\cal
A}(\det P^u )$.
Let $k\in{\cal A}(K_X)$ be the Chern connection of the canonical bundle, and let
$\lambda:=a \otimes k^{\vee}$ be the induced connection in $L:=\det E$.  We
denote by
${\cal L}:=(L,\bar\partial_\lambda)$ the holomorphic structure defined by
$\lambda$. Now
identify the affine space ${\cal A}(\bar\delta(P^u))$ with the space ${\cal
A}_{\lambda\otimes k^{\otimes 2}}(E\otimes K_X)$ of connections in
$E\otimes K_X$
which induce $\lambda\otimes k^{\otimes 2}=a\otimes k$ in $\det(E\otimes
K_X)$, and
identify $A^0(\Sigma^+(P^u))$ with $A^0(E\otimes K_X)\oplus A^0(E)=A^0(E \otimes
K_X)\oplus A^{02}(E\otimes K_X)$.
\vspace{3mm}
\begin{pr} Fix an integrable connection $a\in{\cal A}(\det E\otimes K_X)$.
A pair
$(A,\varphi+\alpha)\in {\cal A}_{\lambda\otimes k^{\otimes 2}}(E\otimes
K_X)\times
[A^0(E \otimes K_X)\oplus A^{02}(E\otimes K_X)]$ solves the
$PU(2)$-monopole equations
$(SW^\sigma_a)$ if and only if $A$ is integrable and one of the following
conditions is
satisfied:
$$\begin{array}{ccccc}I)\  \alpha=0,& \bar\partial_A\varphi=0& and& i\Lambda_g
F_A^0+\frac{1}{2}(\varphi\bar\varphi)_0=0\\
 II)\  \varphi=0,&  \partial_A\alpha=0& and &i\Lambda_g
F_A^0-\frac{1}{2}*(\alpha\wedge\bar\alpha)_0=0\ .
\end{array}$$
\end{pr}
\vspace{3mm}

Note that solutions $(A,\varphi)$ of type  I   give rise to holomorphic
pairs $({\cal
F}_A,\varphi)$, consisting of a holomorphic  structure in $F:=E\otimes K$ and a
holomorphic section  $\varphi$ in ${\cal F}_A$. The remaining equation
$i\Lambda_g
F_A^0+\frac{1}{2}(\varphi\bar\varphi)_0=0$  can  again be interpreted as the
vanishing condition for a moment map for the ${\cal G}_0$-action in the
space of pairs
$(A,\varphi)\in {\cal A}_{\lambda\otimes k^{\otimes 2}}(F)\times A^0(F)$.
We shall
study the corresponding stability condition in the next section.

The analysis of the solutions of type II  can be reduced  to the
investigation of
the type  I  solutions: Indeed, if $\varphi=0$ and $\alpha\in
A^{02}(E\otimes K_X)$
satisfies $\partial_A\alpha=0$, we see that the section
$\psi:=\bar\alpha\in A^0(\bar E)$ must be holomorphic, i.e. it satisfies
$\bar\partial_{A\otimes  [a^{\vee}]}\psi=0$.  On the other hand one has
$-*(\alpha\wedge\bar\alpha)_0=*(\bar\alpha\wedge
\bar{\bar\alpha})_0=(\psi\bar\psi)_0$.

\subsection{Vortex equations and stable oriented pairs}

Let $(X,g)$ be a compact K\"ahler manifold of arbitrary dimension, and let
$E$ be a
differentiable vector bundle of rank $r$, endowed with a fixed holomorphic
structure ${\cal
L}:=(L,\bar\partial_{\cal L})$ in $L:=\det E$.

An \ub{oriented} \ub{pair} of type $(E,{\cal L})$ is a pair $({\cal
E},\varphi)$, consisting of a
holomorphic structure ${\cal E}=(E,\bar\partial_{\cal E})$ in $E$ with
$\bar\partial_{\det{\cal
E}}=\bar\partial_{\cal L}$, and a holomorphic section $\varphi\in H^0({\cal
E})$. Two oriented
pairs are isomorphic if they are equivalent under the natural action of the
group $SL(E)$ of
differentiable automorphisms of $E$ with determinant 1.

An oriented pair $({\cal E},\varphi)$ is \ub{simple} if its stabilizer in
$SL(E)$ is contained in
the center $\Z_r\cdot\id_E$ of $SL(E)$; it is \ub{strongly} \ub{simple}  if
this stabilizer is
trivial.
\begin{pr} {}[OT5] There exists a (possibly non-Hausdorff) complex analytic
orbifold ${\cal
M}^{si}(E,{\cal L})$ parameterizing isomorphism classes of simple oriented
pairs of type
$(E,{\cal L})$. The open subset  ${\cal M}^{ssi}(E,{\cal L})\subset{\cal
M}^{si}(E,{\cal L})$ of
classes of strongly simple pairs is a complex analytic space, and the
points ${\cal
M}^{si}(E,{\cal L})\setminus {\cal M}^{ssi}(E,{\cal L})$ have neighborhoods
modeled on
$\Z_r$-quotients.
\end{pr}
Now fix a Hermitian background metric $H$ in $E$. In this section we use
the symbol
($SU(E)$) $U(E)$ for the groups  of (special) unitary    automorphisms of
$(E,H)$, and not
for the  bundles of (special) unitary automorphisms.

Let $\lambda$ be the Chern connection
associated with the Hermitian holomorphic bundle $({\cal L},\det H)$. We
denote by $\bar{\cal
A}_{\bar\partial_\lambda}(E)$ the affine space of semiconnections in $E$
which induce the
semiconnection $\bar\partial_\lambda=\bar\partial_{\cal L}$ in  $L=\det E$,
and we write
${\cal A}_\lambda(E)$ for the space of  unitary connections in $(E,H)$
which induce $\lambda$
in $L$. The map $A\longmapsto \bar\partial_A$ yields an identification ${\cal
A}_\lambda(E)\map \bar{\cal A}_{\bar\partial_\lambda}(E)$, which endows the
affine space
${\cal A}_{\lambda}(E)$ with a complex structure. Using this identification
and the Hermitian
metric $H$, the product ${\cal A}_{\lambda}(E)\times A^0(E)$ becomes --
after suitable
Sobolev completions -- an infinite dimensional K\"ahler manifold.  The map
$$m:{\cal A}_\lambda(E)\times A^0(E)\map A^0(su(E))
$$
defined by $m(A,\varphi):=\Lambda_g
F_A^0-\frac{i}{2}(\varphi\bar\varphi)_0$ is a
\ub{moment} \ub{map} for the $SU(E)$-action on the K\"ahler manifold ${\cal
A}_\lambda(E)
\times A^0(E)$.

We denote by ${\cal H}_\lambda(E):=\{(A,\varphi)\in{\cal
A}_\lambda(E)\times A^0(E)|\
F_A^{02}=0,\ \bar\partial_A\varphi=0\}$ the space of integrable pairs, and
by  ${\cal
H}_\lambda(E)^{si}$ the open subspace of    pairs $(A,\varphi)\in{\cal
H}_\lambda(E)$ with
$(\bar\partial_A,\varphi)$ simple.   The quotient
$${\cal V}_\lambda(E):=\qmod{{\cal H}_\lambda(E)\cap m^{-1}(0)}{SU(E)}\ ,\ \
(\ {\cal V}^*_\lambda(E):=\qmod{{\cal H}^{si}_\lambda(E)\cap
m^{-1}(0)}{SU(E)}\ )
$$
is called the moduli space  of (irreducible) \ub{projective} \ub{vortices}.
Note that a vortex
$(A,\varphi)$ is irreducible iff $SL(E)_{(A,\varphi)}\subset \Z_r\id_E$.
Using again an
infinite dimensional version of the theory of symplectic quotients   (as in
the abelian
case),  one gets   a homeomorphism
 $$
j:{\cal V}_\lambda(E)\textmap{\simeq} \qmod{{\cal H}_\lambda^{ps}(E)}{SL(E)}
$$
where ${\cal H}_\lambda^{ps}(E)$ is the   subspace of ${\cal H}_\lambda(E)$
consisting of pairs whose $SL(E)$-orbit meets the vanishing locus of the
moment map.  ${\cal
H}_\lambda^{ps}(E)$ is in general not open, but ${\cal H}_\lambda^{s}(E):={\cal
H}_\lambda^{ps}(E)\cap {\cal H}^{si}_\lambda(E)$  is open, and
restricting   $j$ to  ${\cal
V}^*_\lambda(E)$ yields an isomorphism of real analytic orbifolds
$${\cal V}^*_\lambda(E)\textmap{\simeq} \qmod{{\cal
H}_\lambda^{s}(E)}{SL(E)}\subset {\cal
M}^{si}(E,{\cal L}) \ .
$$

 The image
$${\cal M}^{s}(E,{\cal L}):= \qmod{{\cal H}_\lambda^{s}(E)}{SL(E)}$$
of this isomorphism can be identified  with the  set of isomorphism classes
of simple
oriented holomorphic pairs $({\cal E},\varphi)$ of type $(E,{\cal L})$,
with the property  that
${\cal E}$ admits a Hermitian metric with $\det h=\det H$  which solves the
\ub{projective}
\ub{vortex} \ub{equation}
$$i\Lambda_g F_h^0+\frac{1}{2}(\varphi\bar\varphi^h)_0=0\ .
$$
Here   $F_h$ is the curvature of the Chern connection of $({\cal E},h)$.\\

The set ${\cal M}^{s}({\cal E},{\cal L})$ has a purely holomorphic
description as the
subspace of elements $[{\cal E},\varphi]\in{\cal M}^{si}(E,{\cal L})$ which
satisfy a suitable
\ub{stability} condition.

This condition is rather complicated for bundles $E$ of rank $r>2$, but it
becomes very simple
when $r=2$.

Recall that, for any torsion free  coherent sheaf ${\cal F}\ne 0$ over a
$n$-dimensional
K\"ahler  manifold $(X,g)$, one defines the $g$-\ub{slope} of ${\cal F}$ by
$$\mu_g({\cal
F}):=\frac{c_1(\det {\cal F} )\cup [\omega_g]^{n-1}}{\rk({\cal F})}\ .$$

A holomorphic bundle ${\cal E}$ over $(X,g)$ is called
\ub{slope}-\ub{stable} if
$\mu_g({\cal F})<\mu_g({\cal E})$ for all proper coherent  subsheaves
${\cal F}\subset
{\cal E}$.  The bundle
${\cal E}$ is \ub{slope}-\ub{polystable} if it decomposes as a direct sum
${\cal E}=\oplus{\cal
E}_i$ of slope-stable bundles with $\mu_g({\cal E}_i)=\mu_g({\cal E})$.
\begin{dt} Let $({\cal E},\varphi)$ be an oriented pair of type $({\cal
E},{\cal L})$ with $\rk E=2$
over a K\"ahler manifold $(X,g)$. The pair $({\cal E},\varphi)$ is
\ub{stable} if $\varphi=0$ and
${\cal E}$ is slope-stable, or $\varphi\ne 0$ and the divisorial component
$D_\varphi$ of the
zero-locus $Z(\varphi)\subset X$ satisfies $\mu_g({\cal O}_X(D))<\mu_g(E)$.
The pair $({\cal
E},\varphi)$ is \ub{polystable} if it is stable or $\varphi=0$ and ${\cal
E}$ is slope-polystable.
\end{dt}
{\bf Example:} Let $D\subset X$ be an effective divisor defined by a
section $\varphi\in
H^0({\cal O}_X(D))\setminus\{0\}$, and put ${\cal E}:={\cal O}_X(D)\oplus
[{\cal L}\otimes{\cal
O}_X(-D)]$. The pair $({\cal E},\varphi)$ is stable iff $\mu_g({\cal
O}_X(2D))<\mu_g({\cal
L})$.\\

The following result gives a metric characterization of polystable oriented
pairs.

\begin{thry} {}[OT5] Let $E$ be a differentiable vector bundle of rank 2
over $(X,g)$ endowed with
a Hermitian holomorphic structure $({\cal L},l)$ in $\det E$. An oriented
pair of type $(E,{\cal
L})$ is polystable iff ${\cal E}$ admits a Hermitian metric $h$ with $\det
h=l$ which solves the
projective vortex equation
$$i\Lambda_g F_h^0+\frac{1}{2}(\varphi\bar\varphi^h)_0=0\ .$$
If $({\cal E},\varphi)$ is stable, then the metric  $h$ is unique.
\end{thry}
This result identifies the subspace ${\cal M}^{s}(E,{\cal L})\subset {\cal
M}^{si}(E,{\cal L})$ as
the subspace of isomorphism classes of \ub{stable} \ub{oriented} \ub{pairs}.

Theorem 3.2.3 can be used to
show that the moduli spaces ${\cal M}^\sigma_{a}$ of $PU(2)$-monopoles on a
K\"ahler surface
have a natural complex geometric description when the connection $a$ is
integrable.
Recall from section 3.1. that in this case ${\cal M}^\sigma_a$ decomposes
as the union of two
Zariski-closed subspaces
$${\cal M}^\sigma_a=({\cal M}^\sigma_a)_I\cup({\cal M}^\sigma_a)_{II}
$$
according to the two conditions I, II in Proposition 3.1.4.   By this
proposition, both terms of
this union can be identified with   moduli spaces of projective vortices.
Using again the symbol
$^*$ to denote  subsets of points with central stabilizers, one gets the
following
Kobayashi-Hitchin type description of $({\cal M}^\sigma_a)^*$ in terms of
stable oriented
pairs.
\begin{thry}{}[OT5], [T2] If $a\in{\cal A}(\det P^u )$ is integrable, the
moduli space  ${\cal
M}^\sigma_a$ decomposes as an union ${\cal M}^\sigma_a=({\cal
M}^\sigma_a)_I\cup({\cal
M}^\sigma_a)_{II}$ of two Zariski closed subspaces isomorphic with moduli
spaces of
projective vortices, which intersect along the Donaldson moduli space
${\cal D}(\bar\delta(P^u))$. There are natural real analytic isomorphisms
$$({\cal M}^\sigma_a)_I^*\textmap{\simeq} {\cal M}^{s}(E\otimes K_X,{\cal
L}\otimes{\cal
K}_X^{\otimes 2})\ ,\ \
({\cal M}^\sigma_a)_{II}^*\textmap{\simeq} {\cal M}^{s}(E^{\vee},{\cal
L}^{\vee}) \ ,
$$
where ${\cal L}$ denotes the holomorphic structure in $\det E=\det P^u
\otimes K_X^{\vee}$
defined by $\bar\partial_a$ and the canonical  holomorphic structure in $K_X$.
\end{thry}
\vspace{8mm}
{\bf Example:}  (R. Plantiko) On $\P^2$, endowed with the standard
Fubini-Study metric $g$,
we consider the $Spin^{U(2)}(4)$-structure $\sigma:P^u\map P_g$ defined by
the standard
$Spin^c(4)$-structure  $\tau:P^c\map  P_g$ and the $U(2)$-bundle $E$ with
$c_1(E)=7$,
$c_2(E)=13$, and we fix an integrable connection $a\in{\cal A}(\det P^u)$. This
$Spin^{U(2)}(4)$-structure is characterized by $c_1(\det(P^u))=4$,
$p_1(\bar\delta(P^u))=-3$, and the  bundle $F:=E\otimes K_{\P^2}$ has Chern
classes $c_1(F)=1$, $c_2(F)=1$. It is easy to see that every stable
oriented pair $({\cal
F},\varphi)$ of type $(F,{\cal O}_{\P^2}(1))$ with $\varphi\ne 0$ fits into
an exact  sequence
of the form
$$
 0\map {\cal O} \textmap{\varphi} {\cal F}\map J_{Z(\varphi)} \otimes {\cal
O}_{\P^2}(1)\
\map 0\ ,$$
where ${\cal F}={\cal T}_{\P^2}(-1)$  and the zero locus $Z(\varphi)$ of
$\varphi$ consists of
a simple point $z_\varphi\in\P^2$. Two such pairs $({\cal F},\varphi)$,
$({\cal F},\varphi')$
define the same point in the moduli space ${\cal M}^s(F,{\cal
O}_{\P^2}(1))$ if and only if
$\varphi'=\pm \varphi$. The resulting identification
$${\cal M}^s(F,{\cal O}_{\P^2}(1))=
\qmod{H^0({\cal T}_{\P^2}(-1))}{\{\pm \id\}}$$
is a complex analytic isomorphism.

Since every polystable pair of type $(F,{\cal O}_{\P^2}(1))$ is actually
stable, and
since there are no polystable oriented pairs of type $(E^{\vee},{\cal
O}_{\P^2}(-7))$,
Theorem 3.2.4 yields a real analytic isomorphism
$${\cal M}^\sigma_a=\qmod{H^0({\cal T}_{\P^2}(-1))}{\{\pm \id\}}$$
where the origin corresponds to the unique stable oriented pair of the form
$({\cal
T}_{\P^2}(-1),0)$.
 The quotient $\qmod{H^0({\cal
T}_{\P^2}(-1))}{\{\pm \id\}}$ has a natural algebraic compactification
${\cal C}$, given by
the cone over the image of $\P(H^0({\cal T}_{\P^2}(-1)))$ under the
Veronese map to
$\P(S^2H^0({\cal T}_{\P^2}(-1)))$.    This
compactification coincides with the Uhlenbeck compactification $\overline{{\cal
M}^\sigma_a}$ (see section 2.2, [T3]). More precisely, let
$\sigma':P'^u\map  P_g$ be the $Spin^{U(2)}(4)$-structure with $\det
P'^u=\det P^u$ and
$p_1(P'^u)=1$. This structure is associated with $\tau$ and the
$U(2)$-bundle $E'$ with Chern classes $c_1(E')=7$, $c_2(E')=12$. The moduli
space ${\cal
M}^{\sigma'}_a$ consists of one (abelian) point, the class of the abelian
solution corresponding
to the
\ub{stable} oriented pair $({\cal O}_{\P^2}\oplus{\cal
O}_{\P^2}(1),\id_{{\cal O}_{\P^2}})$ of
type $(E'\otimes K_{\P^2},{\cal O}_{\P^2}(1))$. ${\cal M}^{\sigma'}_a$ can be identified with
the moduli space ${\cal W}^\tau_{\frac{i}{2\pi}F_a}$ of
$ [\frac{i}{2\pi}F_a ]$-twisted abelian
\sw monopoles. Under the identification ${\cal C}=\overline{{\cal
M}^\sigma_a}$, the vertex
of the cone corresponds to the unique Donaldson point which is given by the
stable oriented
pair
$({\cal T}_{\P^2}(-1),0)$. The base of the cone corresponds to the space ${\cal
M}^{\sigma'}_a\times
\P^2$ of ideal monopoles concentrated in one point.\\
\\

We want to close this section by explaining the stability concept which
describes the subset
${\cal M}^{s}_X(E,{\cal L})\subset {\cal M}^{si}_X(E,{\cal L})$ in the
general case $r\geq 2$. This
stability concept    does \ub{not} depend on the choice of parameter and
the corresponding
moduli spaces can be interpreted as "master spaces" for holomorphic pairs
(see next section);
  in the projective framework they admit Gieseker type compactifications [OST].

We shall find this stability concept by relating the $SU(E)$-moment map $m:{\cal
A}_\lambda(E)\times A^0(E)\map A^0(su(E))$ to the universal family of
$U(E)$-moment  maps
$m_t:{\cal A}(E)\times A^0(E)\map A^0(u(E))$ defined by
$$m_t(A,\varphi):=\Lambda_g F_A
-\frac{i}{2}(\varphi\bar\varphi)+\frac{i}{2}t\id_E\ ,
$$
where $t\in A^0$ is an arbitrary real valued function. Given $t$, we
consider the following
system of equations
$$\left\{\begin{array}{ccc}
F_A^{02}&=&0\\
\bar\partial_A\varphi&=&0\\
i\Lambda_g F_A +\frac{1}{2}(\varphi\bar\varphi)&=&\frac{t}{2}\id_E
\end{array}\right. \eqno{(V_t)}
$$
for pairs $(A,\varphi)\in A(E)\times A^0(E)$. Put $\rho_t:=\frac{1}{4\pi
n}\int\limits_X t\omega_g^n $. To explain our first result, we have to
recall some classical
stability concepts for holomorphic pairs.\\

For any   holomorphic bundle ${\cal E}$ over $(X,g)$   denote by ${\cal
S}({\cal E})$ the set of
reflexive subsheaves ${\cal F}\subset {\cal E}$ with $0<\rk({\cal
F})<\rk({\cal E})$, and  for
a fixed section $\varphi\in H^0({\cal E})$   put
$${\cal S}_\varphi({\cal E}):=\{{\cal F}\in{\cal S}({\cal E})|\ \varphi\in
H^0({\cal F})\}\ .
$$
Define real numbers  $\underline{m}_g({\cal E})$ and $\overline{m}_g({\cal
E},\varphi)$ by
$$\underline{m}_g({\cal E}):=\max(\mu_g({\cal E}),\sup\limits_{{\cal
F}'\in{\cal S}({\cal E})}
\mu_g({\cal F}'))\ ,\
\overline{m}_g({\cal E},\varphi):=\inf\limits_{{\cal F}\in{\cal S}_\varphi({\cal
E})} \mu_g(\qmod{{\cal E}}{{\cal F}})\ .
$$
A bundle ${\cal E}$ is $\varphi$-\ub{stable} in the sense of  S. Bradlow when
$\underline{m}_g({\cal E})<\overline{m}_g({\cal E},\varphi)$.
Let $\rho\in\R$ be any real parameter.  A holomorphic pair $({\cal
E},\varphi)$ is called
$\rho$-\ub{stable}  if $\rho$ satisfies  the inequality
$$\underline{m}_g({\cal E})<\rho<\overline{m}_g({\cal E},\varphi) \ .
$$
The pair $({\cal E},\varphi)$ is  $\rho$-\ub{polystable} if  it is
$\rho$-stable or ${\cal
E}$-splits holomorphically  as ${\cal E}={\cal E}'\oplus{\cal E}''$ such
that $\varphi\in
H^0({\cal E}')$, $({\cal E}',\varphi)$ is $\rho$-stable and ${\cal E}''$
is a slope-polystable
vector bundle with $\mu_g({\cal E}'')=\rho$ [Br].  Let $GL(E)$ be the group
of bundle
automorphisms of $E$. With these definitions one proves [OT1]  (see  [Br]
for the case of a
constant function
$t$):
\begin{pr} The complex orbit $(A,\varphi)\cdot GL(E)$ of an integrable pair
$(A,\varphi)\in {\cal A}(E)\times A^0(E)$ contains a solution of $(V_t)$ if
and only if
the pair $({\cal E}_A,\varphi)$ is $\rho_t$-polystable.
\end{pr}

Let us now   fix again a Hermitian metric $H$ in $E$ and an integrable
connection $\lambda$ in
the Hermitian line bundle $L:=(\det E,\det H)$. Consider  the system of
equations
$$\left\{\begin{array}{ccc}
F_A^{02}&=&0\\
\bar\partial_A\varphi&=&0\\
i\Lambda_g F_A^0 +\frac{1}{2}(\varphi\bar\varphi)_0&=&0
\end{array}\right. \eqno{(V^0)}
$$
for pairs $(A,\varphi)\in {\cal A}_\lambda(E)\times A^0(E)$. Then one can prove
\begin{pr} Let $(A,\varphi)\in {\cal A}_\lambda(E)\times A^0(E)$ be an
integrable pair. The
following assertions are  equivalent:\\
i) The complex orbit $SL(E)\cdot(A,\varphi)$ contains a solution of $(V^0)$.\\
ii) There exists a function $t\in A^0$ such that the $GL(E)$-orbit
$GL(E)\cdot(A,\varphi)$
contains a solution of $(V_t)$.\\
iii) There exists a real number $\rho$ such that the pair $({\cal
E}_A,\varphi)$ is
$\rho$-polystable.
\end{pr}
\begin{co} The open subspace ${\cal M}^{s}(E,{\cal L})\subset {\cal
M}^{si}(E,{\cal L})$   is the
set of isomorphism classes of simple oriented pairs which are
$\rho$-polystable for some
$\rho\in\R$.
\end{co}

\begin{re} There exist stable oriented pairs $({\cal E},\varphi)$ whose
stabilizer
with respect to the $GL(E)$-action is  of positive dimension. Such
pairs cannot be $\rho$-stable   for any  $\rho\in\R$.
\end{re}

Note that the moduli spaces ${\cal M}^{s}(E,{\cal L})$ have a natural
$\C^*$-action defined by
$z\cdot [{\cal E},\varphi]:=[{\cal E},z^{\frac{1}{r}}\varphi]$. This action
is well defined since
$r$-th roots of unity are contained in the complex gauge group $SL(E)$.

There exists an equivalent definition for stability of
oriented pairs, which does not use  the parameter dependent stability
concepts of [Br]. The
fact that it is expressible in terms of $\rho$-stability is related to the
fact that the
moduli spaces ${\cal M}^{s}(E,{\cal L})$ are master spaces for moduli spaces of
$\rho$-stable pairs.
\subsection{Master spaces and the coupling principle}
Let $X\subset\P^N_\C$ be a smooth  complex projective variety with
hyperplane bundle ${\cal
O}_X(1)$. All degrees and Hilbert polynomials of coherent sheaves will be
computed
corresponding to these data.

We fix a torsion-free sheaf ${\cal E}_0$ and a holomorphic line bundle
${\cal L}_0$ over $X$, and
we choose a Hilbert polynomial $P_0$. By $P_{\cal F}$ we denote the Hilbert
polynomial of a
coherent sheaf ${\cal F}$. Recall that any non-trivial torsion free
coherent sheaf ${\cal F}$
admits a unique subsheaf ${\cal F}_{\max}$  for which $\frac{P_{\cal
F'}}{\rk{\cal F}'}$ is
maximal and whose rank is maximal among  all subsheaves ${\cal F}'$ with
$\frac{P_{\cal
F'}}{\rk{\cal F}'}$ maximal.

An ${\cal L}_0$-\ub{oriented} \ub{pair} of type $(P_0,{\cal E}_0)$ is a
triple $({\cal
E},\varepsilon,\varphi)$ consisting of a torsion free coherent sheaf ${\cal
E}$ with
determinant isomorphic to ${\cal L}_0$ and Hilbert polynomial $P_{\cal
E}=P_0$, a
homomorphism $\varepsilon:\det{\cal E}\map {\cal L}_0$, and a morphism
$\varphi:{\cal
E}\map{\cal E}_0$. The homomorphisms $\varepsilon$ and $\varphi$ will be
called the
\ub{orientation} and the  \ub{framing} of the oriented pair. There is an
obvious equivalence
relation for such pairs. When $\ker\varphi\ne 0$, we set
$$\delta_{{\cal E},\varphi}:=P_{\cal E}-\frac{\rk{\cal
E}}{\rk[\ker(\varphi)_{\max}]}
P_{\ker(\varphi)_{\max}}\ .
$$

An oriented pair $({\cal E},\varepsilon,\varphi)$ is   \ub{semistable} if
one of the
following conditions is  satisfied:\\
1. $\varphi$ is injective.\\
2. $\varepsilon$ is an isomorphism, $\ker \varphi\ne 0$,  $\delta_{{\cal
E},\varphi}\geq
0$, and for all non-trivial subsheaves ${\cal F}\subset {\cal E}$ the
following inequality
holds
$$\frac{P_{\cal F}}{\rk{\cal F}}-\frac{\delta_{{\cal E},\varphi}}{\rk{\cal
F}}\leq
\frac{P_{\cal E}}{\rk{\cal E}} -\frac{\delta_{{\cal E},\varphi}}{\rk{\cal E}}\ .
$$

The corresponding \ub{stability} concept is slightly more complicated
[OST].  Note that the
(semi)stability definition above does not depend on a parameter. It is,
however, possible to
express (semi)stability
 in terms of the parameter   dependent Gieseker-type stability  concepts of
[HL2].  E.g.,
$({\cal E},\varepsilon,\varphi)$ is semistable iff $\varphi$ is injective,
or ${\cal E}$ is
Gieseker semistable, or there exists a rational polynomial $\delta$ of
degree smaller than
$\dim X$  with
positive leading coefficient, such that $({\cal E},\varphi)$ is
$\delta$-semistable in the sense
of [HL2].

For all stability  concepts introduced so far there exist analogous notions of
slope-(semi)stability. In  the
special case when the reference sheaf ${\cal E}_0$ is the trivial sheaf ${\cal O}_X$, slope
stability is the algebro-geometric analog of the stability concept
associated with the
projective vortex equation.
\begin{thry} {}[OST] There exists a \ub{projective} scheme ${\cal
M}^{ss}(P_0,{\cal E}_0,{\cal
L}_0)$ whose closed points correspond to  {\rm gr}-equivalence classes of
Gieseker semistable
${\cal L}_0$-oriented pairs of type $(P_0,{\cal E}_0)$. This scheme
contains an open subscheme
${\cal M}^{s}(P_0,{\cal E}_0,{\cal L}_0)$ which is a coarse moduli space
for stable ${\cal
L}_0$-oriented pairs.
\end{thry}
 It is also possible to construct moduli spaces for stable oriented pairs
where the orienting
line bundle is allowed to vary [OST]. This generalization is important in
connection with
Gromov-Witten invariants for Grassmannians [BDW].

Note that ${\cal M}^{ss}(P_0,{\cal E}_0,{\cal L}_0)$ possesses a natural
$\C^*$-action, given
by
$$z\cdot [{\cal E},\varepsilon,\varphi]:=[{\cal E},\varepsilon,z\varphi]\ ,
$$
whose fixed point set can be explicitly described.  The fixed point locus
$[{\cal
M}^{ss}(P_0,{\cal E}_0,{\cal L}_0)]^{\C^*}$ contains two distinguished
subspaces, ${\cal
M}_{0}$ defined by the equation $\varphi=0$, and ${\cal M}_{\infty}$ defined by
$\varepsilon=0$. ${\cal M}_{0}$ can be identified with the Gieseker scheme
${\cal
M}^{ss}(P,{\cal L}_0)$ of equivalence classes of semistable ${\cal
L}_0$-oriented torsion free
coherent sheaves with Hilbert polynomial $P_0$. The subspace ${\cal
M}_{\infty}$    is the
Grothendieck Quot-scheme $Quot_{P_{{\cal E}_0}-P_0}^{{\cal E}_0,{\cal
L}_0}$ of quotients of
${\cal E}_0$ with fixed determinant isomorphic with $(\det {\cal
E}_0)\otimes {\cal
L}_0^{\vee}$ and Hilbert polynomial  $P_{{\cal E}_0}-P_0$.

In the terminology  of   [BS],  ${\cal
M}_0$ is  the source
${\cal M}_{source}$ of the $\C^*$-space ${\cal M}^{ss}(P_0,{\cal E}_0,{\cal
L}_0)$, and
${\cal M}_{\infty}$ is its sink when non-empty.

The remaining subspace of the fixed point locus
$${\cal M}_R:=[{\cal M}^{ss}(P_0,{\cal
E}_0,{\cal L}_0)]^{\C^*}\setminus[{\cal M}_0\cup{\cal M}_{\infty}]\ ,$$
 the so-called space of
\ub{reductions},  consists of objects which are  of the same type but
essentially of lower
rank.

Note that  the
Quot scheme ${\cal M}_{\infty}$ is empty if $\rk({\cal E}_0)$ is smaller
than the rank $r$
of the sheaves ${\cal E}$ under consideration, in which case the sink of
the moduli space is
a closed subset of the space of reductions.

Recall from [BS] that
 the closure of a general $\C^*$-orbit connects a point in ${\cal
M}_{source}$ with a point in
${\cal M}_{sink}$, whereas closures of special orbits connect points of
other parts of the fixed
point set.

 The flow generated by the $\C^*$-action can therefore be used to relate
data associated
with
${\cal M}_{0}$ to data associated with ${\cal M}_{\infty}$ and ${\cal M}_R$.

The technique of computing  data on ${\cal M}_{0}$ in terms of  ${\cal
M}_R$ and ${\cal
M}_{\infty}$ is a very general principle which we call \ub{coupling}
\ub{and} \ub{reduction}.
This principle has already been  described in a gauge theoretic framework
in section 2.2 for
relating monopoles and instantons. However, the essential ideas may
probably be best
understood   in an abstract  Geometric Invariant Theory  setting, where one
has a very
simple and clear picture.
\\

Let $G$ be a complex reductive group, and consider a linear  representation
$\rho_A:G\map
GL(A)$ in a finite dimensional vector space  $A$. The induced action
$\bar\rho_A \map{\rm
Aut}(\P(A))$ comes with a natural linearization in ${\cal O}_{\P(A)}(1)$,
hence we have a
stability concept, and thus  we can form the GIT quotient
$${\cal M}_0:=\P(A)^{ss}//G\ .$$
Suppose we want to compute "correlation functions"
$$\Phi_I:=\langle \mu_I,[{\cal M}_0]\rangle \ ,
$$
i.e. we want to evaluate suitable products of canonically defined
cohomology classes $\mu_i$
on the fundamental class $[{\cal M}_0]$ of  ${\cal M}_0$. Usually the
$\mu_i$'s are slant
products  of characteristic classes of a "universal  bundle"
$\raisebox{0.25ex}{$\scriptstyle|$}\hskip -1.5mm{{\cal E}}_0$  on ${\cal M}_0\times X$ with
homology classes of $X$ . Here
$X$ is a compact manifold,  and
$\raisebox{0.29ex}{${\scriptstyle |}$}\hskip -1.5mm{{\cal E}}_0$    comes from a
 tautological bundle  $\raisebox{0.29ex}{${\scriptstyle |}$}\hskip
-1.5mm{\tilde{\cal E}}_0$ on $A\times X$ by applying Kempf's Descend Lemma.

The main idea is now to couple the original problem with a simpler  one,
and to use the
$\C^*$-action which occurs naturally in the resulting GIT quotients to
express the original
correlation functions in terms of simpler data. More precisely, consider another
representation
$\rho_B:G \map GL(B)$ with GIT quotient ${\cal M}_{\infty}:=\P(B)^{ss}//G$.
The direct sum
$\rho:=\rho_A\oplus\rho_B$ defines a naturally linearized $G$-action on the
projective
space
$\P(A\oplus B)$. We call the corresponding quotient
$${\cal M}:=\P(A\oplus B)^{ss}//G
$$
the \ub{master} \ub{space} associated with the coupling of $\rho_A$ to $\rho_B$.

 The space ${\cal M}$ comes with a natural $\C^*$-action, given by
$$z\cdot [a,b]:=[a,z\cdot b]\ ,
$$
and the union ${\cal M}_0\cup{\cal M}_{\infty}$ is a closed subspace of the
fixed point locus
${\cal M}^{\C^*}$.

Now make the simplifying assumptions that ${\cal M}$ is smooth and
connected, the
$\C^*$-action is free outside ${\cal M}^{\C^*}$, and suppose that the
cohomology classes
$\mu_i$ extend to ${\cal M}$. This condition is always satisfied if the
$\mu_i$'s were
obtained by the procedure   described above, and if Kempf's lemma applies
to the pull-back
bundle
$p_A^*(\raisebox{0.29ex}{${\scriptstyle |}$}\hskip -1.5mm{\tilde{\cal
E}}_0)$  and
provides a bundle on ${\cal M}\times X$  extending
$\raisebox{0.29ex}{${\scriptstyle
|}$}\hskip -1.5mm{{\cal E}}_0$.

Under these assumptions, the complement
$$ {\cal M}_R:={\cal M}^{\C^*}\setminus ({\cal M}_0\cup{\cal M}_{\infty})
$$
is a closed submanifold of ${\cal M}$, disjoint from ${\cal M}_0$, and
${\cal M}_{\infty}$. We
call ${\cal M}_R$ the manifold of \ub{reductions} of the master space.  Now
remove a
sufficiently small $S^1$-invariant tubular neighborhood $U$ of ${\cal
M}^{\C^*}\subset {\cal
M}$, and consider the $S^1$-quotient $W:=\qmod{[{\cal M}\setminus
U]}{S^1}$. This is a
compact manifold whose boundary is the union of the projectivized normal bundles
$\P(N_{{\cal M}_0})$ and $\P(N_{{\cal M}_{\infty}})$, and a differentiable
projective
fiber space $P_R$ over ${\cal M}_R$.  Note that in general  $P_R$ has  no
natural holomorphic
structure. Let $n_0$, $n_\infty$ be the complex dimensions of the fibers of
$\P(N_{{\cal
M}_0})$, $\P(N_{{\cal M}_{\infty}})$, and let
$u\in H^2(W,\Z)$ be the first Chern class of the $S^1$-bundle dual   to
${\cal M}\setminus
U\map W$. Let $\mu_I$ be a class as above. Then, taking into account
orientations, we compute:
$$\Phi_I:=\langle \mu_I,[{\cal M}_0]\rangle=\langle \mu_I\cup
u^{n_0},[\P(N_{{\cal
M}_0})]\rangle=\langle \mu_I\cup u^{n_0},[\P(N_{{\cal
M}_\infty})]\rangle -\langle \mu_I\cup u^{n_0}, [P_R]\rangle \ .
$$

In this way the coupling principle reduces the calculation of the original
correlation
functions on
${\cal M}_0$  to computations on ${\cal M}_\infty$ and on the manifold of
reductions ${\cal
M}_R$. A particular important case occurs when  the GIT problem  given by
$\rho_B$ is
trivial, i.e. when $\P(B)^{ss}=\emptyset$. Under these circumstances the
functions $\Phi_I$
are completely determined by data associated with the manifold of
reductions ${\cal M}_R$.

Of course, in realistic situations, our simplifying assumptions are seldom
satisfied, so that
one has to modify the basic idea in a suitable way.

One of the realistic situations which we have in mind is the coupling of
coherent sheaves
with morphisms into a fixed reference sheaf ${\cal E}_0$.      In this
case, the original problem is the classification of stable torsion-free
sheaves, and the
corresponding Gieseker scheme ${\cal M}^{ss}(P_0,{\cal L}_0)$ of ${\cal
L}_0$-oriented
semistable sheaves of Hilbert polynomial $P_0$ plays the role of the
quotient ${\cal M}_0$.
The corresponding master spaces are the moduli spaces
${\cal M}^{ss}(P_0,{\cal E}_0,{\cal L}_0)$ of semistable ${\cal
L}_0$-oriented pairs of type
$(P_0,{\cal E}_0)$.

Coupling with ${\cal E}_0$-valued homomorphisms $\varphi:{\cal E}\map {\cal
E}_0$ leads to
two essentially different situations, depending on the rank $r$  of the sheaves
${\cal E}$ under consideration:

1. When $\rk({\cal E}_0)<r$, the framings $\varphi:{\cal E}\map{\cal E}_0$
can never
be injective, i.e. there are no semistable homomorphisms. This case
correspond  to the GIT
situation ${\cal M}_\infty=\emptyset$.

2. As soon as $\rk({\cal E}_0)\geq  r$, the framings $\varphi$ can become
injective, and the
Grothendieck schemes  $Quot_{P_{{\cal E}_0}-P_0}^{{\cal E}_0,{\cal L}_0}$
appear in the
master space ${\cal M}^{ss}(P_0,{\cal E}_0,{\cal L}_0)$. These Quot schemes
are the analoga
of  the quotients ${\cal M}_{\infty}$ in the  GIT situation.

In both cases the spaces of reductions are  moduli spaces of   objects
which are of the same
type but essentially of lower rank.

Everything can be made very explicit when the base manifold is a curve $X$
with a trivial
reference sheaf ${\cal E}_0={\cal O}_X^{\oplus k}$. In the case $k< r$, the
master spaces relate
correlation functions of moduli spaces of semistable bundles with fixed
determinant to data
associated with reductions. When $r=2$, $k=1$, the manifold of reductions
are symmetric
powers of the base curve, and the coupling principle can be used to prove
the \ub{Verlinde}
\ub{formula}, or to compute the volume and the characteristic numbers (in
the smooth case)  of
the moduli spaces of  semistable bundles.

The general case $k\geq r$ leads to a method for the computation   of
\ub{Gromov}-\ub{Witten}
invariants for Grassmannians.  These
invariants can be regarded as correlation functions of suitable Quot
schemes [BDW], and the
coupling principle relates them to data associated with reductions and
moduli spaces of
semistable bundles. In this case one needs a master space ${\cal
M}^{ss}(P_0,{\cal
E}_0,{\cal L})$ associated with a Poincar\'e line bundle ${\cal L}$ on
$\Pic(X)\times X$ which
 set theoretically is the union over ${\cal L}_0\in\Pic(X)$ of the master
spaces   ${\cal
M}^{ss}(P_0,{\cal E}_0,{\cal L}_0)$ [OST]. One could try to  prove the
\ub{Vafa}-\ub{Intriligator} formula along these  lines.

Note that the use of master spaces allows us to avoid the sometimes messy
investigation
of   \ub{chains} of \ub{flips}, which occur whenever one considers the family of
all possible
$\C^*$-quotients of the master space  [T], [BDV].

The coupling principle has been   applied in two further situations.

Using the coupling of vector bundles with twisted endomorphisms, A. Schmitt
has recently
constructed projective moduli spaces of \ub{Hitchin} \ub{pairs} [S]. In the
case of curves
and twisting with the canonical bundle, his master spaces are natural
compactifications of
the moduli spaces introduced in [H].

Last but not least, the coupling principle can  also be used in certain
gauge theoretic
situations:

The coupling of instantons on 4-manifolds with Dirac-harmonic spinors  has
been described
in detail in chapter 2. In this case the instanton moduli spaces are the
original moduli
spaces ${\cal M}_{0}$, the Donaldson polynomials are the original correlation
functions to compute, and the moduli spaces of
$PU(2)$-monopoles    are master
spaces for the coupling with   spinors. One is again in the special
situation where ${\cal
M}_{\infty}=\emptyset$, and the manifold of reductions is a union of moduli
spaces of
twisted abelian monopoles. In order to compute the contributions  of the
abelian moduli
spaces to the correlation functions, one  has to give explicit
descriptions of the master
space in an
$S^1$-invariant neighborhood of the   abelian locus.

Finally consider again the Lie group $G=Sp(n)\cdot S^1$ and the
$PSp(n)$-monopole
equations
$(SW^\sigma_a)$ for a   $Spin^{Sp(n)\cdot S^1}(4)$-structure
$\sigma:P^G\map  P_g$ in $(X,g)$ and an abelian connection $a$ in the
associated $S^1$-bundle
(see Remark 2.1.2).   Regarding the compactification of the moduli space  ${\cal M}^\sigma_a$
as master space associated with the coupling of $PSp(n)$-instantons to harmonic
spinors, one should get a relation between Donaldson $PSp(n)$-theory and
Seiberg-Witten
type theories.

\newpage
\section{\bf References}
\vspace{10 mm}
\parindent 0 cm

[AHS] Atiyah M.; Hitchin N. J.; Singer I. M.: {\it Selfduality in
four-dimensional Riemannian geometry}, Proc. R. Soc. Lond. (A) 362, 425-461
(1978)

[BPV] Barth, W.; Peters, C.; Van de Ven, A.: {\it Compact complex surfaces},
Springer Verlag (1984)

[B] Bauer, S.: {\it Some non-reduced moduli of bundles and Donaldson
invariants for
Dolgachev surfaces}, J. Reine Angew. Math. 442, 149-161 (1996)

[BDW]  Bertram, A.;   Daskalopoulos, G.;   Wentworth, R.: {\it
Gromov invariants for holomorphic maps from Riemann surfaces to
Grassmannians},  J. Amer.
Math. Soc.  9 ( 2),  529-571  (1996)

[BS] Bia\l ynicki-Birula, A.;  Sommese A.: {\it
Quotients by $\C^*$ and $SL_2(\C)$ actions}
Trans.  Am.  Math.  Soc.  279, 773-800  (1983)

[Br] Bradlow, S. B.: {\it Special metrics and stability for holomorphic
bundles with global sections}, J. Diff. Geom. 33, 169-214 (1991)

[BrDW] Bradlow, S.;   Daskalopoulos, G.;   Wentworth, R.: {\it
Birational equivalences of vortex moduli}, Topology  35 (3), 731-748  (1996)

[Bru] Brussee, R.: {\it Some $ {\cal C}^{\infty}$-properties of K\"ahler
surfaces}, preprint,
alg-geom 9503004 (1995)

[BW]  Bryan, J. A.; Wentworth, R.: {\it The multi-monopole equations for
K\"ahler surfaces},
Turkish Journal of Math.   20 (1), 119-128  (1996)

[DH]   Dolgachev, I.;  Hu, Y.:  {\it Variation of geometric invariant
theory quotients},
preprint

 [D] Donaldson, S.: {\it The Seiberg-Witten equations and 4-manifold
topology}, Bull. AMS
  33 (1),  45-70  (1996)

[DK] Donaldson, S.; Kronheimer, P.: {\it The Geometry of
four-manifolds}, Oxford Science Publications   (1990)

[FS] Fintushel, R.; Stern, R.:  {\it Rational blowdowns of smooth
4-manifolds},  preprint,
alg-geom/9505018  (1995)

[FU] Freed D. S.;  Uhlenbeck, K.:
{\it Instantons and Four-Manifolds },
Springer-Verlag (1984)

[FM] Friedman, R.; Morgan, J. W.: {\it Algebraic surfaces and \sw
invariants}, preprint,
alg-geom/9509007  (1995)

[F] Furuta, M.: {\it Monopole equation and the $\frac{11}{8}$-conjecture},
preprint,
RIMS, Kyoto (1995)

[G]  Gieseker, D.: {\it On the moduli of vector bundles on an algebraic
surface},
Ann.  Math.   106, 45-60  (1977)

[GS] Guillemin, V.; Sternberg, S.: {\it Birational equivalence in the
symplectic category},  Inv. math. 97, 485-522 (1989)

[H1] Hitchin, N.: {\it  Harmonic spinors}, Adv. in Math. 14, 1-55 (1974)

[H2]   Hitchin, N.: {\it The self-duality equations on a Riemann surface},
Proc. London  Math. Soc.     55 (3), 59-126  (1987)

[HKLR] Hitchin, N.; Karlhede, A.; Lindstr\"om, U.; Ro\v cek, M.: {\it
Hyperk\"ahler metrics and
supersymmetry}, Commun.  Math.  Phys.  108, 535-589 (1987)

[HH] Hirzebruch, F.; Hopf,  H.: {\it Felder von Fl\"achenelementen in
4-dimensiona\-len  Mannigfaltigkeiten}, Math. Ann. 136, 156-172 (1958)

[HL]    Huybrechts, D.;   Lehn, M.: {\it Framed modules and their moduli},
Internat.\ J.\ Math.   6, 297-324  (1995)

 [JPW] Jost, J.; Peng, X.;  Wang, G.: {\it Variational aspects of the
Seiberg-Witten functional},
Calc. Var. Partial Differential Equations 4 ( 3), 205-218 (1996)

[K] Kobayashi, S.: {\it Differential geometry of complex vector bundles},
Princeton University Press  (1987)

[KN]  Kobayashi, S.;  Nomizu, K.:
{\it Foundations of differential geometry I,} John Wiley \& Sons (1963)

[KM1] Kronheimer, P.; Mrowka, T.: {\it  Recurrence relations and asymptotics for
four-manifold invariants}, Bull. AMS 30 (2), 205-221 (1994)

[KM2] Kronheimer, P.; Mrowka, T.: {\it The genus of embedded surfaces in
the projective plane}, Math. Res. Lett. 1, 797-808 (1994)

[LM] Labastida, J. M. F.;  Marino, M.: {\it Non-abelian monopoles on
 four manifolds}, preprint,
Departamento de Fisica de Particulas,   Santiago de Compostela, April
 (1995)

[LMi] Lawson, H. B. Jr.; Michelson, M. L.: {\it Spin Geometry}, Princeton
University Press (1989)

[L1]  Le Brun, C.: {\it On the scalar curvature of complex surfaces},
Geometric and Functional
Analysis   5 (3), 619-628 (1995)

[L2]  Le Brun, C.: {\it Einstein metrics and Mostow rigidity},
Math. Res. Lett. 2,  1-8   (1995)

[LL] Li, T.; Liu, A.: {\it General wall crossing formula}, Math. Res. Lett.
2,
  797-810  (1995)

[LT] L\"ubke, M.; Teleman, A.: {\it The Kobayashi-Hitchin
correspondence},
 World Scientific Publishing Co.  (1995)

[Lu] Lupa\c scu, P.:{\it Riemannian metrics on surfaces of K\"ahler type},
preprint,
Z\"urich (1996)

[M]  Miyajima, K.: {\it Kuranishi families of
vector bundles and algebraic description of
the moduli space of Einstein-Hermitian
connections},   Publ. R.I.M.S. Kyoto Univ.  25,
  301-320 (1989)

[Mo] Morgan, J.: {\it The \sw equations and applications to the topology of
smooth
four-manifolds}, Princeton University Press  (1996)

[MFK] Mumford, D.; Fogarty, J.; Kirwan, F.: {\it Geometric invariant
theory},  Springer Verlag
(1994)

[OST]  Okonek, Ch.; Schmitt, A.; Teleman, A.: {\it Master spaces for stable
pairs}, preprint,
alg-geom/9607015 (1996)

[OT1]   Okonek, Ch.; Teleman, A.: {\it The Coupled Seiberg-Witten
Equations, Vortices, and Moduli Spaces of Stable Pairs},    Int. J. Math.
  6 (6), 893-910 (1995)

[OT2] Okonek, Ch.; Teleman, A.: {\it Les invariants de Seiberg-Witten
et la conjecture de Van De  Ven}, Comptes Rendus Acad. Sci. Paris, t.
321,  S\'erie I, 457-461 (1995)

[OT3] Okonek, Ch.; Teleman, A.: {\it Seiberg-Witten invariants and
rationality of complex surfaces}, Math. Z., to appear

[OT4] Okonek, Ch.; Teleman, A.: {\it Quaternionic monopoles},  Comptes
Rendus Acad. Sci. Paris, t. 321, S\'erie I, 601-606 (1995)

[OT5]  Okonek, Ch.;  Teleman, A.: {\it Quaternionic monopoles},
Commun.\ Math.\ Phys.  180 (2), 363-388  (1996)

[OT6]  Okonek, Ch..;  Teleman, A.: {\it Seiberg-Witten invariants for
manifolds with $b_+=1$,
and the universal wall crossing formula},
Int. J. Math. 7 (6), 811-832 (1996)

[OV]  Okonek, Ch.;  Van de Ven, A.: {\it $\Gamma$-type invariants associated to
$PU(2)$-bundles and the differentiable structure of  Barlow's surface},
Invent. Math. 95,
601-614 (1989)

[PT1]  Pidstrigach, V.; Tyurin, A.: {\it Invariants of the smooth
structure of an algebraic surface arising from the Dirac operator},
Russian Acad. Izv. Math.    40 (2), 267-351 (1993)

[PT2]  Pidstrigach, V.; Tyurin, A.: {\it Localisation of the Donaldson
invariants along the
Seiberg-Witten classes}, Russian Acad. Izv., to appear

[S] Schmitt, A.: {\it Projective moduli for Hitchin pairs}, preprint,
Z\"urich  (1996)

[Ta1] Taubes, C.: {\it The Seiberg-Witten invariants and symplectic forms},
Math.
Res. Lett. 1, 809-822  (1994)

[Ta2] Taubes, C.: {\it The Seiberg-Witten invariants and the Gromov
invariants}, Math.
Res. Lett. 2, 221-238  (1995)

[Ta3] Taubes, C.: {\it $SW\Rightarrow Gr$. From the Seiberg-Witten equations to
pseudo-holomorphic curves}, J. Amer. Math. Soc. 9 (3), 845-918  (1996)

[T1] Teleman, A.: {\it Non-abelian Seiberg-Witten theory},
Habilitationsschrift, Universit\"at
Z\"urich (1996)

[T2] Teleman, A.: {\it Non-abelian Seiberg-Witten theory and  stable
oriented  pairs}, Int. J.
Math., to appear

[T3] Teleman, A.: {\it Moduli spaces of $PU(2)$-monopoles},
preprint, Universit\"at Z\"urich  (1996)

[Te] Teleman, A. M.: {\it Seiberg-Witten functionals}, in preparation

[Th]  Thaddeus,  M.: {\it Geometric invariant theory and flips}, preprint.

[W] Witten, E.: {\it Monopoles and four-manifolds}, Math.  Res.
Lett. 1,  769-796  (1994)
\vspace{1cm}\\
Authors'  addresses : \\

 Institut f\"ur Mathematik, Universit\"at Z\"urich,  Winterthurerstr. 190, \\
CH-8057 Z\"urich,  \\
{  e-mail}: okonek@math.unizh.ch\ \ \ \  teleman@math.unizh.ch

\end{document}